# Diffusion models with physics-guided inference for solving partial differential equations

Bing Yi[1*], Jia Liu[1], Jinyang Fu[2,3*], Xiang Peng[4]

[1] School of Traffic and Transportation Engineering, Central South University, Changsha, China

[2] School of Civil Engineering, Central South University, Changsha, China

[3] National Engineering Research Center of High-speed Railway Construction Technology, Changsha, China

[4] College of Mechanical Engineering, Zhejiang University of Technology, Hangzhou, China

**Abstract:** Diffusion models have recently emerged as powerful stochastic frameworks for high-dimensional inference and generation. However, existing applications to partial differential equations (PDEs) predominantly rely on physics-informed training strategies, which tightly couple learning with specific governing equations and limit generalization across problem settings. In this work, we propose a diffusion model with physics-guided inference for solving PDEs, in which the diffusion model is trained using standard data-driven procedures, while physical laws are incorporated exclusively during the reverse inference stage. The reverse diffusion dynamics is guided by a PDE residual energy function, combined with Gaussian smoothing and explicit boundary enforcement, yielding a physically consistent stochastic iteration that is independent of the training process. From a numerical standpoint, the proposed framework can be interpreted as a diffusion-inspired implicit solver that converges to the PDE solution even when initialized from random noise and perturbed by stochastic fluctuations. The method is validated on classical PDE equation such as Poisson, Diffusion, and Burgers equations with varying coefficients. Numerical results demonstrate robust convergence, high accuracy, and strong generalization without retraining, highlighting the proposed framework as a unified alternative to classical numerical solvers and physics-informed neural networks.

## 1. Introduction

Partial differential equations (PDEs) form the mathematical backbone of modeling and simulation in engineering and applied sciences. Fundamental physical processes such as diffusion, transport, wave propagation, and conservation laws are commonly expressed through PDEs, including Poisson, diffusion, and Burgers-type equations. Accurate and efficient numerical solution of PDEs is essential in a wide range of applications, such as heat transfer, fluid dynamics, solid mechanics, electrochemical systems, and multi-physics engineering design [1–3].

Over the past several decades, extensive efforts have been devoted to developing reliable and efficient methods for solving PDEs. Despite the diversity of approaches, these methods share a common goal: obtaining accurate solutions under complex geometries, boundary conditions, material parameters, and operating regimes. As engineering systems become increasingly high-dimensional, multiscale, and data-rich, the inherent limitations of traditional numerical paradigms—particularly in handling the curse of dimensionality and multi-physics coupling—have become more pronounced. This exigency necessitates a transition toward integrated mathematical frameworks that synergize classical first principles with data-driven intelligence, driving a new wave of methodological innovation [4].

### 1.1 Classical numerical methods: accuracy and robustness

Classical numerical methods, including finite difference methods (FDM), finite element methods

(FEM), finite volume methods (FVM), and spectral methods, have long served as the foundation of PDE solving in computational mechanics and engineering analysis [1,2,5–8]. These methods are derived from rigorous mathematical principles and provide strong guarantees in terms of consistency, stability, and convergence. Boundary conditions and governing equations are explicitly enforced, making such methods highly reliable for industrial and safety-critical applications.

However, the accuracy and robustness of classical solvers come at the cost of significant computational expense. High-resolution meshes, small time steps, and repeated linear or nonlinear system solves are often required, especially for parametric studies, optimization loops, uncertainty quantification, and real-time control [9,10]. Although advanced techniques such as multigrid solvers and domain decomposition methods have substantially improved efficiency [11,12], the computational burden remains prohibitive in many practical scenarios.

## 1.2 Data-driven learning methods: efficiency with data dependence

With the rapid growth of computational power and the availability of large-scale simulation and experimental datasets, data-driven learning methods have emerged as an alternative paradigm for accelerating PDE solving. Early neural network based surrogate models approximate PDE solutions from labeled data, enabling fast inference once trained. More recently, operator learning frameworks have attracted considerable attention by learning mappings between infinite-dimensional function spaces rather than discrete solution snapshots.

Representative approaches include Deep Operator Networks (DeepONet) [13] and neural operator architectures such as Fourier Neural Operators (FNO) and graph-based neural operators [14,15]. These methods have demonstrated impressive performance in learning parametric PDE families and generalizing across discretization. By shifting computational effort from online inference to offline training, operator learning significantly reduces solution time in many-query settings.

Nevertheless, data-driven methods exhibit several limitations. Their performance strongly depends on the availability of large, high-quality training datasets, which are often generated using expensive numerical solvers. Moreover, generalization may deteriorate when extrapolating to unseen boundary conditions, coefficients, or forcing terms, limiting robustness in real-world applications [4].

## 1.3 Physics-informed learning: embedding equations into training

To address the lack of physical consistency in purely data-driven models, physics-informed neural networks (PINNs) were introduced by incorporating PDE residuals and boundary conditions directly into the training loss function [16,17]. PINNs eliminate the need for labeled data and provide a unified framework for forward and inverse problems.

Numerous extensions have been proposed to improve PINN performance, including conservative formulations [18], adaptive loss weighting strategies [19], and spatiotemporal domain decomposition strategies for large-scale or long-term problems [20]. Despite these advances, PINNs often suffer from slow convergence, optimization stiffness, and sensitivity to hyperparameters. These issues become particularly restrictive when resolving sharp gradients, convection-dominated regimes, or long-term temporal dynamics [21]. Furthermore, a major limitation of the standard PINN framework is the lack of adaptability: any alteration in physical parameters or boundary conditions typically necessitates a complete retraining of the model.

## 1.4 Generative diffusion models: strong generalization with limited physics

Generative diffusion models have recently emerged as a powerful class of stochastic generative models for high-dimensional data [22,23]. By formulating generation as a gradual denoising process governed by stochastic differential equations, diffusion models exhibit stable training behavior and strong generalization capabilities.

These properties have motivated growing interest in applying diffusion models to scientific computing tasks, including uncertainty quantification, inverse problems, and advancements in the mathematical foundations of score-based models [24–26]. Compared to supervised learning approaches, diffusion models are less sensitive to limited or noisy data and can capture complex solution distributions.

However, diffusion models are inherently data-driven. Incorporating strict physical laws—such as differential operators and boundary conditions—into the generative process remains challenging. Existing guidance strategies, such as classifier-based or score-based guidance, rely on auxiliary learned models and typically enforce physical constraints only implicitly through data, rather than explicitly through governing equations [27]. As a result, physical consistency during inference is not guaranteed.

### 1.5 Motivation and contributions

The above discussion reveals a clear methodological gap. Classical numerical solvers provide high accuracy and physical rigor but are computationally expensive. Data-driven and physics-informed learning methods improve efficiency but often require extensive retraining and careful optimization. Diffusion models offer strong generalization and robustness but lack principled mechanisms for incorporating physical laws during inference.

To bridge this gap, we propose a diffusion model with physics-guided inference for solving partial differential equations. In contrast to physics-informed training strategies, the diffusion model in this work is trained using standard data-driven procedures. Physical laws are incorporated exclusively during the reverse inference process through a PDE energy function, Gaussian smoothing, and explicit boundary enforcement.

The main contributions of this work are summarized as follows:

1) A physics-guided inference framework that decouples diffusion model training from PDE constraints;

2) A stochastic variational interpretation linking diffusion models, Gibbs measures, and PDE solutions;

3) Robust PDE solving under varying coefficients without retraining;

4) A unified alternative to classical numerical solvers and physics-informed neural networks.

## 2. Physics-Guided Diffusion Inference for PDEs

### 2.1 Problem formulation and objective

We consider partial differential equations (PDEs) defined on a bounded spatial domain $\Omega \subset \mathbb{R}^d$ with boundary $\partial\Omega$:

$$\mathcal{L}u + \mathcal{N}(u) = f \ \text{ in } \Omega \tag{1}$$

subject to boundary conditions:

$$\mathcal{B}u = g \ \text{ on } \partial\Omega \tag{2}$$

where $\mathcal{L}$ denotes a linear differential operator, $\mathcal{N}(\cdot)$ represents nonlinear terms, and $f$ is a given source.

Classical numerical solvers approximate Eq. (1) and (2) through deterministic discretization and iterative schemes, while learning-based approaches typically aim to learn the solution operator or enforce the PDE residual during training.

In contrast, the goal of this work is to solve PDEs through a stochastic inference process, where:

- the diffusion model is trained in a purely data-driven manner,
- physical laws are incorporated only during inference,
- convergence to the PDE solution is achieved through physics-guided denoising.

## 2.2 Standard diffusion model training

The denoising diffusion probabilistic model (DDPM) framework is adopted. Let $u_0 \sim p_{\text{data}}(u)$ denote samples of PDE solutions obtained from high-fidelity numerical solvers.

**Forward diffusion:** The forward process gradually corrupts the data by injecting Gaussian noise over T steps. The state at arbitrary time step $t$ can be expressed in a closed form directly from the clean data $u_0$:

$$u_t = \sqrt{\bar{\alpha}_t} u_0 + \sqrt{1 - \bar{\alpha}_t} \xi_t, \;\; \xi_t \sim \mathcal{N}(0, I) \tag{3}$$

where $\{\alpha_t\}_{t=1}^{T}$ is a predefined noise schedule, and $\bar{\alpha}_t = \prod_{i=1}^{t} \alpha_i$ is the cumulative noise magnitude such that $\overline{\alpha_T} \approx 0$.

**Training objective:** A neural network $\epsilon_\theta(u_t, t)$, implemented using a U-Net architecture, is trained to predict the noise $\epsilon$ added to $u_t$. The training objective is to minimize the mean squared error:

$$\mathcal{L}(\theta) = \mathbb{E}_{u_0, t, \epsilon}[||\epsilon - \epsilon_\theta(u_t, t)||^2] \tag{4}$$

Importantly, no physical constraints or boundary conditions are imposed during training. Consequently, the model learns a purely data-driven representation of the solution space, maintaining strict independence from specific PDE formulations at the training stage.

## 2.3 Reverse diffusion and probabilistic inference

In the DDPM framework, sampling corresponds to simulating a reverse-time stochastic process that progressively removes noise.

In discrete form, the reverse update is given by:

$$u_{t-1} = \mu_\theta(u_t, t) + \sigma_t \xi_t, \xi_t \sim \mathcal{N}(0, I) \tag{5}$$

where mean $\mu_\theta(\cdot)$ is parameterized using the trained noise predictor $\epsilon_\theta$.

In the continuous-time limit, this process is elegantly described by the reverse-time stochastic differential equation (SDE) [22, 28]:

$$du_t = \left[-\frac{1}{2}\beta(t) u_t - \beta(t) \nabla_u \log p_t(u_t)\right] dt + \sqrt{\beta(t)} d\bar{W}_t \tag{6}$$

where $\beta(t)$ is the noise schedule and $\bar{W}_t$ is a reverse-time Wiener process. The continuous and discrete schedules are linked via $\alpha_t = exp(-\int_0^t \beta(s) ds)$, ensuring the process remains variance preserving.

The score function $\nabla_u \log p_t(u_t)$ points toward high-density regions of the data distribution—is approximated by the trained network as:

$$\nabla_u \log p_t(u_t) \approx -\frac{\epsilon_\theta(u_t, t)}{\sqrt{1 - \bar{\alpha}_t}} \tag{7}$$

This relationship allows us to substitute the neural network's output directly into the SDE to perform sampling.

## 2.4 Classifier-guided diffusion

To condition generation on external information, classifier-guided diffusion augments the reverse drift by an additional score term [25,27-29]:

$$du_t = \left[-\frac{1}{2}\beta(t)u_t - \beta(t)\nabla_u \log p_t(u_t) - \lambda\nabla_u \log p(y|u_t)\right]dt + \sqrt{\beta(t)}d\bar{W}_t \tag{8}$$

where $p(y|u)$ is a classifier likelihood (or any differentiable constraint) and $\lambda > 0$ controls guidance strength.

This formulation establishes a critical principle: the inference trajectory of a diffusion model can be dynamically steered by augmenting the score function with the gradient of an external energy or log-density. This mechanism provides the formal basis for our physics-guided approach, where physical laws are enforced during sampling rather than training.

## 2.5 PDE energy function

To incorporate physical laws during inference, we define a PDE energy functional $E(u)$ that characterizes the residual of the governing equations:

$$E(u) = \frac{1}{2}\int_{\Omega} (\mathcal{L}u + \mathcal{N}(u) - f)^2 \, \mathrm{d}\Omega \tag{9}$$

where $E(u) = 0$ if and only if $u$ satisfies the PDE under the prescribed boundary conditions.

The variational gradient $\nabla_u E(u)$ represents a deterministic descent direction in the functional space, effectively driving the generated samples toward the manifold of physically admissible solutions.

## 2.6 Physics-guided reverse diffusion

Inspired by classifier guidance, we augment the reverse SDE by replacing the task-specific classifier score with the PDE energy gradient. The resulting physics-guided reverse SDE is:

$$\mathrm{d}u_t = \left[-\frac{1}{2}\beta(t)u_t - \beta(t)\nabla_u \log p_t(u_t) - \lambda\nabla_u E(u_t)\right]\mathrm{d}t + \sqrt{\beta(t)}\,\mathrm{d}\bar{W}_t \tag{10}$$

where:

- The first two terms perform data-driven denoising, leveraging the learned prior.
- The third term enforces physical consistency via the energy gradient.
- $\lambda$ is a guidance scale that balances the generative prior with physical constraints.

As $t \to 0$, the noise schedule $\beta(t) \to 0$ and the influence of the learned prior diminishes, while the PDE-guided term dominates, ensuring convergence toward a physically valid solution.

## 2.7 Gaussian smoothing and boundary enforcement

Applying $\nabla_u E(u)$ directly to noisy intermediate states $u_t$ can lead to numerical instability. To ensure a smooth and meaningful guidance field, we introduce a Gaussian smoothing operator $\mathcal{G}_\sigma$:

$$\nabla_u E(u) \longrightarrow \nabla_u E(\mathcal{G}_\sigma u) \tag{11}$$

where $\mathcal{G}_\sigma$ filters out high-frequency noise that would otherwise corrupt the physical gradients. Furthermore, boundary conditions are explicitly enforced after each update step through a projection operator $\mathcal{P}_\mathcal{B}$:

$$u \leftarrow \mathcal{P}_\mathcal{B}(u) \tag{12}$$

It includes the Dirichlet, Neumann, or periodic constraints. Together, these steps maintain the structural integrity and physical consistency of the solution throughout the stochastic sampling process.

## 2.8 Energy-based interpretation and convergence

As the process approaches convergence ($t \approx 0$), the dynamics are governed primarily by the energy gradient and residual noise. The system simplifies to a Langevin-type SDE:

$$du = -\nabla_u E(u)dt + \sqrt{2\varepsilon}\, dW_t \tag{13}$$

where $\varepsilon$ represents the terminal stochasticity. The associated Fokker–Planck equation admits a stationary distribution known as the Gibbs measure:

$$p_\infty(u) \propto exp\left(-\frac{E(u)}{\varepsilon}\right) \tag{14}$$

This result provides a rigorous statistical foundation for our method: the inference process converges in distribution to a state centered at the PDE solution. In the limit $\varepsilon \to 0$, the distribution collapses to a Dirac delta mass concentrated on the exact solution of the PDE, demonstrating the theoretical robustness of the physics-guided inference.

# 3. PDE-Specific Forms and Inference Equations

In this section, we instantiate the proposed physics-guided diffusion inference framework for representative classes of partial differential equations, including elliptic, parabolic, and nonlinear hyperbolic–parabolic problems. For each PDE, we explicitly construct the corresponding energy functional and derive the associated inference dynamics used during the reverse diffusion process.

## 3.1 Poisson equation

We first consider the Poisson equation with diffusion coefficient, subject to prescribed Dirichlet or Neumann boundary conditions on $\partial\Omega$, the governing equation is formulated as:

$$-\nabla \cdot \left(\kappa \nabla u\right) = f \text{ in } \Omega \tag{15}$$

where $u$ denotes the solution field and $f$ represents the source term. In this setting, $\kappa$ serves as a constant scalar coefficient governing the global diffusion rate.

The governing equation is enforced through a residual-based energy function:

$$E_P(u) = \frac{1}{2}\| - \nabla \cdot (\kappa \nabla u) - f\|^2 \tag{16}$$

This formulation penalizes violations of the Poisson operator in an $L_2$ sense and admits the exact PDE solution as its unique minimizer under appropriate boundary conditions. Compared with variational formulations based on primal energy minimization, the residual-based energy is particularly convenient for diffusion inference, as it directly yields a gradient flow in the solution space without introducing additional test functions.

During reverse diffusion, the solution field is updated according to:

$$\mathrm{d}u_t = \left[-\frac{1}{2}\beta(t)u_t - \beta(t)\nabla_u \log p_t(u_t) - \lambda \nabla_u E_P(\mathcal{G}_\sigma u_t)\right]\mathrm{d}t + \sqrt{\beta(t)}\,\mathrm{d}\bar{W}_t \tag{17}$$

Here, $\mathcal{G}_\sigma$ denotes Gaussian smoothing, which regularizes high-frequency noise and ensures that differential operators are applied to sufficiently smooth fields. The third term corresponds to a gradient descent step on the Poisson energy, while the stochastic perturbation enables exploration of the solution landscape during inference. Crucially, the boundary conditions are enforced as hard constraints through the projection operator $\mathcal{P}_\mathcal{B}$ after each update, ensuring the state is consistently mapped back into the admissible space. From a numerical perspective, this update can be interpreted as a stochastic solver that converges to the Poisson solution even when initialized from random noise.

## 3.2 Heat diffusion equation

To demonstrate the framework's capability in handling time-dependent processes, we first consider the classical transient diffusion equation:

$$\frac{\partial u}{\partial t} - \nabla \cdot \left(\alpha \nabla u\right) = f(x,t) \tag{18}$$

where $u(\mathrm{x}, \mathrm{t})$ represents the temperature field, $\alpha$ denotes the thermal diffusivity, and $f(x,t)$ is an external source term. Standard numerical methods typically solve Eq. (18) by marching through discrete time steps, which can lead to cumulative temporal errors.

**Space-Time Stationary Transformation:** We reframe the temporal dimension $t$ as a spatial-like coordinate $y \in [0, \mathrm{Y}]$, thereby transforming the transient evolution into a stationary boundary value problem on a unified space-time domain $D = \Omega \times [0, \mathrm{Y}]$:

$$u_y - \alpha u_{xx} = f(x,y) \tag{19}$$

In this formulation, the physical constraints are partitioned into the governing equation within $D$ and specific boundary requirements on $\partial D$. The initial condition is treated as a Dirichlet constraint at the boundary with $y = 0$:

$$u(x,0) = h(x) \tag{20}$$

while the spatial boundary conditions are defined along the boundaries of $x \in \partial\Omega$:

$$u(x,y)|_{x\in\partial\Omega} = g(x,y) \tag{21}$$

**Residual-based Energy Construction:** The governing physics within the domain are then enforced through a residual-based energy function:

$$E_H(u) = \frac{1}{2}\|u_y - \alpha u_{xx} - f\|^2 \tag{22}$$

This function serves as the objective for physical guidance. Unlike soft-penalty methods, the initial and boundary conditions $h(x)$ and $g(x,y)$ are explicitly enforced through the projection operator $\mathcal{P}_{\mathcal{B}}$ introduced in Section 2.7, ensuring that every inference step remains strictly within the admissible solution space.

The diffusion framework acts as a physics-guided solver that optimizes the entire space-time field $u(x,y)$ simultaneously. To distinguish from the physical time $t$, we denote the iterative inference stage by virtual time $\tau$. The physics-guided reverse SDE is expressed as:

$$du_\tau = \left[-\frac{1}{2}\beta(\tau)u_\tau - \beta(\tau)\nabla_u \log p_\tau(u_\tau) - \lambda\nabla_u E_H(\mathcal{G}_\sigma u_\tau)\right]d\tau + \sqrt{\beta(\tau)}\, d\bar{W}_\tau \tag{23}$$

Followed by the hard-constraint projection $u_\tau = \mathcal{P}_{\mathcal{B}}(u_\tau)$. Here, $\mathcal{G}_\sigma$ denotes Gaussian smoothing, which regularizes high-frequency noise and ensures that the differential operators are applied to sufficiently smooth fields. Even for linear diffusion, this smoothing is essential during the early stages of inference when the state $u_\tau$ is dominated by stochastic perturbations. By treating the entire space-time evolution as a unified manifold, the method by passes step-by-step error accumulation. The solution is obtained through deterministic inference driven by the energy gradient, which is integrated with diffusion-inspired optimization to provide enhanced robustness compared with explicit schemes while avoiding iterative nonlinear solvers.

## 3.3 Burgers equation

Finally, we consider the nonlinear Burgers' equation, a fundamental model for capturing the interaction between nonlinear advection and viscous diffusion:

$$\frac{\partial u}{\partial t} + u \cdot \nabla u - \nu\Delta u = 0, \qquad x \in \Omega, t \in [0, \mathrm{T}] \tag{24}$$

where $\nu$ represents the kinematic viscosity. This equation is characterized by the formation of sharp gradients and shocks, posing a significant challenge for traditional numerical solvers.

**Space-Time Stationary Transformation:** Consistent with our holistic approach, we transform the transient Burgers' equation into a stationary boundary value problem by mapping the temporal dimension t to a spatial-like coordinate $y \in [0, \mathrm{Y}]$:

$$u_y + uu_x - \nu u_{xx} = 0 \tag{25}$$

The physical constraints are defined by the initial state $u(x,0) = h(x)$ and the spatial boundary conditions $u(x,y)|_{x\in\partial\Omega} = g(x,y)$. Similar to the heat equation, these constraints are enforced explicitly via the projection operator $\mathcal{P}_{\mathcal{B}}$ during the inference process.

**Nonlinear Residual-based Energy:** The governing physics within the domain are then enforced through the following residual-based energy function:

$$E_B(u) = \frac{1}{2}\|u_y + uu_x - \nu u_{xx}\|^2 \tag{26}$$

Similarly, the initial and boundary conditions $h(x)$ and $g(x,y)$ are explicitly enforced through the projection operator $\mathcal{P}_{\mathcal{B}}$, ensuring that every inference step remains strictly within the admissible solution space. This energy captures the competition between advective steepening and diffusive smoothing. The variational gradient $\nabla_u E_B$ introduces strong nonlinearity into the inference dynamics, steering the generated field toward a physically valid solution that respects the conservation laws of the system.

The inference process is performed by simulating the physics-guided reverse SDE over the iterative stage $\tau$.

$$du_\tau = \left[-\frac{1}{2}\beta(\tau)u_\tau - \beta(\tau)\nabla_u \log p_\tau(u_\tau) - \lambda\nabla_u E_B(\mathcal{G}_\sigma u_\tau)\right]d\tau + \sqrt{\beta(\tau)}\, d\bar{W}_\tau \tag{27}$$

Followed by the hard-constraint projection $u_\tau = \mathcal{P}_{\mathcal{B}}(u_\tau)$. In this context, the Gaussian smoothing operator $\mathcal{G}_\sigma$ is particularly essential; it stabilizes the guidance field by filtering high-frequency noise that could otherwise lead to numerical divergence near sharp gradients or emerging shocks. Despite the inherent non-convexity of the energy landscape in nonlinear problems, the proposed framework remains stable and convergent. By treating the entire space-time evolution as a unified consistent manifold, the model successfully recovers complex nonlinear dynamics without any specific retraining of the underlying diffusion prior, demonstrating the robust zero-shot generalization of our physics-guided paradigm.

## 3.4 Algorithm summary

The proposed physics-guided diffusion inference procedure for solving partial differential equations is summarized in Algorithm 1. In the practical implementation, the energy gradient is replaced by the product of the physical equation residual (as defined in Eq. (28)) and a physical guidance step size to guide the diffusion model's inference. Crucially, when the residual becomes zero, the physical energy error also vanishes, indicating that the model has converged to the exact solution. We have demonstrated the validity of this simplification through extensive experiments.

$$\mathcal{R}(u) = \mathcal{L}u + \mathcal{N}(u) - f \tag{28}$$

For transient problems, we treat the temporal dimension as a spatial-like coordinate, thereby transforming a one-dimensional transient evolution into a two-dimensional stationary field for unified inference. For steady-state problems, the process initiates from a two-dimensional Gaussian random field, consistent with the terminal distribution of standard diffusion models, and iteratively generates the solution by reversing the diffusion process.

At each reverse step, a pretrained diffusion model first predicts a denoised field, providing a data-driven prior that captures the global solution manifold. This prediction is subsequently refined through physics-guided correction, where a Gaussian smoothing operator suppresses high-frequency noise and a PDE energy gradient enforces consistency with the governing equations. Finally, boundary conditions are explicitly imposed at every iteration via a projection operator to guarantee strict physical admissibility.

Unlike conventional diffusion models that rely on learned classifiers or conditional embeddings for guidance, the proposed method incorporates deterministic PDE operators directly into the inference dynamics. This framework ensures that the stochastic denoising process is continuously projected toward the PDE-constrained solution space. In the zero-noise limit, the algorithm recovers a classical PDE solver, while for finite noise levels it performs stochastic inference that converges in distribution to the Gibbs measure induced by the physical energy.

**Algorithm 1: Physics-Guided Diffusion Inference for PDEs**

Input: Pretrained diffusion model $\epsilon_\theta$, PDE energy functional $E$, guidance time step size $\Delta t$, smoothing parameter $\sigma$, Noise's variance schedule $\beta_t$, $\alpha_t = 1 - \beta_t$, Standard Gaussian noise $\xi_t$, $\bar{\alpha}_t = \prod_{i=1}^{t} \alpha_i$

1. Initialize $u_T \sim \mathcal{N}(0, I)$
2. For $t = T, \ldots, 1$ do:
   - Predict denoised field using the diffusion model (Data-driven prior):

$$u \leftarrow \frac{1}{\sqrt{\alpha_t}}\left(u - \frac{\beta_t}{\sqrt{1-\bar{\alpha}_t}}\epsilon_\theta(\mathrm{x}_t, t)\right) + \sqrt{\beta_t} \cdot \xi_t$$

   - Apply Gaussian smoothing to the intermediate estimate: $u \leftarrow \mathcal{G}_\sigma(u)$.
   - Update using PDE energy gradient to enforce physical consistency:

$$\mathcal{R}(u) \leftarrow \nabla_u E(u) \qquad u \leftarrow u - \Delta t * \mathcal{R}(u)$$

   - Enforce boundary conditions explicitly via projection: $u \leftarrow \mathcal{P}_{\mathcal{B}}(u)$.

   End For
3. Output final PDE solution $u_0$.

# 4. Numerical Experiments

This section evaluates the proposed diffusion model with physics-guided inference on representative partial differential equations, specifically focusing on forward problems with varying coefficients to verify the model's generalization capability. The proposed method is compared with classical numerical solvers, physics-informed neural networks (PINNs), and unguided diffusion models in terms of convergence behavior, accuracy, and generalization capability.

## 4.1 Dataset generation and Experimental Setup

The physical domain is defined as a unit square $\Omega = [0,1]^2$, discretized into a $64 \times 64$ grid to align with the network's resolution requirements. Training datasets were generated using the MATLAB's PDE solver, which employs the Finite Difference Method (FDM) or Finite Element Method (FEM) for high-fidelity integration. Unless otherwise specified, each dataset consists of 4,000 spatiotemporal snapshots, with equation parameters uniformly sampled from defined intervals to ensure comprehensive coverage of the physical regimes. To stabilize the training dynamics and ensure numerical consistency, all solution fields were normalized to $[-1,1]$ via Global Max-Absolute Scaling.

All computations were performed within the MATLAB environment on a workstation equipped with an Intel Core i9-13900HX CPU and an NVIDIA RTX 4060 GPU.

## 4.2 Network Architecture and Training Strategy

We employ a modified U-Net as the backbone for the Denoising Diffusion Probabilistic Model

(DDPM). The architecture is designed for continuous field regression, featuring a 3-level encoder-decoder structure with 64 base channels and residual connections. This setup enables the network to effectively capture multiscale spatial correlations within the solution manifold, mapping a Gaussian noise field $u_T$ to the predicted noise residual $\epsilon_\theta$.

The training process utilizes a linear variance schedule and minimizes the Mean Squared Error (MSE) loss via the Adam optimizer. Importantly, no physical constraints are incorporated during training phase, ensuring the model serves as a generic, data driven prior over solution manifold. The diffusion timesteps and physics-guided step sizes were calibrated based on the specific stiffness and dynamics of each equation. Hyperparameters utilized for network are summarized in Table 1.

Table 1. Hyperparameter configuration of the proposed PDE-guided diffusion model

| Hyperparameter | Value |
|---|---|
| Learning Rate | $10^{-3}$ |
| Batch Size | 64 |
| Noise schedule range ($\beta_t$) | [$10^{-4}$, 0.02] |
| Smoothing parameter ($\sigma$) | 0.9 |

## 4.3 Evaluation Protocol and Baseline Methods

To strictly assess the performance of the proposed physics-guided diffusion framework, we benchmark it against three distinct categories of solvers: high-fidelity classical numerical schemes, physics-informed learning methods, and pure data-driven generative models.

### 4.3.1 Baseline methods

Classical Numerical Solvers: High-resolution numerical solutions serve as the ground truth for error quantification. For one-dimensional problems, we employ the Finite Difference Method (FDM) with a sufficiently fine spatial-temporal grid to ensure numerical stability and accuracy. For two-dimensional domains, solutions are computed using the Finite Element Method (FEM) on a refined mesh to minimize discretization errors and handle complex domain geometries.

Physics-informed learning: We utilize Physics-Informed Neural Networks (PINNs) as the representative baseline for physics-informed learning method. Unlike the proposed generative framework, PINNs treat the PDE solution as a function approximation problem, requiring specific retraining for each new configuration of coefficients or boundary conditions. The network parameters $\theta$ are optimized by minimizing a composite loss function $\mathcal{L}_\theta$ that enforces the PDE governing equations and boundary conditions:

$$\mathcal{L}_\theta = \lambda_{res}\mathcal{L}_{\text{res}} + \lambda_{bc}\mathcal{L}_{\text{bc}} + \lambda_{data}\mathcal{L}_{\text{data}} \tag{29}$$

where $\mathcal{L}_{\text{res}}$ denotes the PDE residual, $\mathcal{L}_{\text{bc}}$ represents the boundary conditions, and $\mathcal{L}_{\text{data}}$ refers to sparse measurements or labeled data points. The coefficients $\lambda$ balance the contribution of each term.

Unguided diffusion model: To isolate the specific contribution of the physics-guided strategy, we evaluate the standard pre-trained diffusion model as an ablation baseline. In this configuration, samples are generated via the reverse diffusion process without any physical guidance or energy-based corrections, relying solely on the learned data distribution.

### 4.3.2 Evaluation metrics

To provide a comprehensive assessment, we employ both quantitative metrics and visual analysis. Global accuracy is quantified via the full-field relative $L_2$ error, as defined in Eq. (30). Visually, we present error heatmaps illustrating the point-wise absolute error $\left(\left|u_{pred} - u_{gt}\right|\right)$ across the domain. Furthermore, we analyze the convergence trajectory to track evolution from initialization to the final

state, and explicitly test generalization capability on coefficient distributions strictly unseen during training stage.

$$\sqrt{\frac{\sum\left(u_{pred}-u_{gt}\right)^2}{\sum u_{gt}^2}} \tag{30}$$

## 4.4 Poisson Equation

To validate the model's capability in solving linear elliptic boundary value problems, we firstly consider the steady-state Poisson equation.

### 4.4.1 Experimental Setup

To rigorously evaluate the accuracy and generalization capability of the proposed method, we constructed a dataset covering diverse range of equation parameters. Training data were generated using a classical numerical solver, with specific coefficient $\kappa$ sampled uniformly at random within the interval $[1.0, 2.0]$. The model was trained on 4,000 training samples for 400 epochs, resulting in 62 iterations per epoch and a total of 24,800 gradient updates based on the batch size in Table 1, employing a decay schedule that halved the learning rate every 100 epochs.

For inference, we designed a robust testing protocol comprising four distinct cases: two used coefficients from within the training interval—but distinct from the training set—to assess interpolation accuracy; the remaining two cases utilized parameters strictly outside the training range to verify extrapolation capability. In these experiments, the diffusion process was fixed at 1000 time steps, while the physics-guided step size was calibrated for each specific case to balance the data prior and physical constraints.

### 4.4.2 Training Dynamics and Forward Diffusion

Prior to evaluating the physics-guided inference solver, we examine the fundamental training phase of the diffusion model. Figure 1 illustrates both the forward diffusion process and the network's training progression. At $t=0$ in Fig 1, the training data consist of smooth solution fields generated via a numerical PDE solver. During training, these two-dimensional fields are progressively corrupted by a forward noise process governed by the Markov chain $q(x_t|x_{t-1})$. As the diffusion step $t$ increases from 0 to 1000, structured physical features are systematically eroded, eventually culminating in an isotropic distribution that approximates pure Gaussian noise.

The objective of the diffusion model is to learn the noise residual at varying diffusion stages, enabling it to subsequently reverse the corruption process. Figure 2 plots the training progress, where the red line tracks the Root Mean Squared Error (RMSE) to provide an intuitive analysis of the convergence behavior. The loss decreases sharply during the initial epochs and reaches a stable plateau after approximately 10,000 iterations. This trend indicates that the U-Net has effectively captured the underlying manifold of the training data distribution, establishing a robust generative prior for the subsequent physics-guided inference.

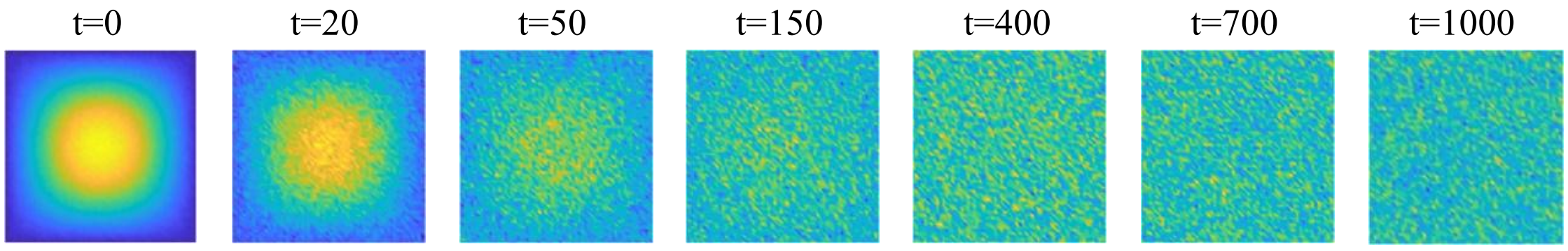

Fig. 1. The forward noising process of the diffusion model for the steady-state Poisson fields

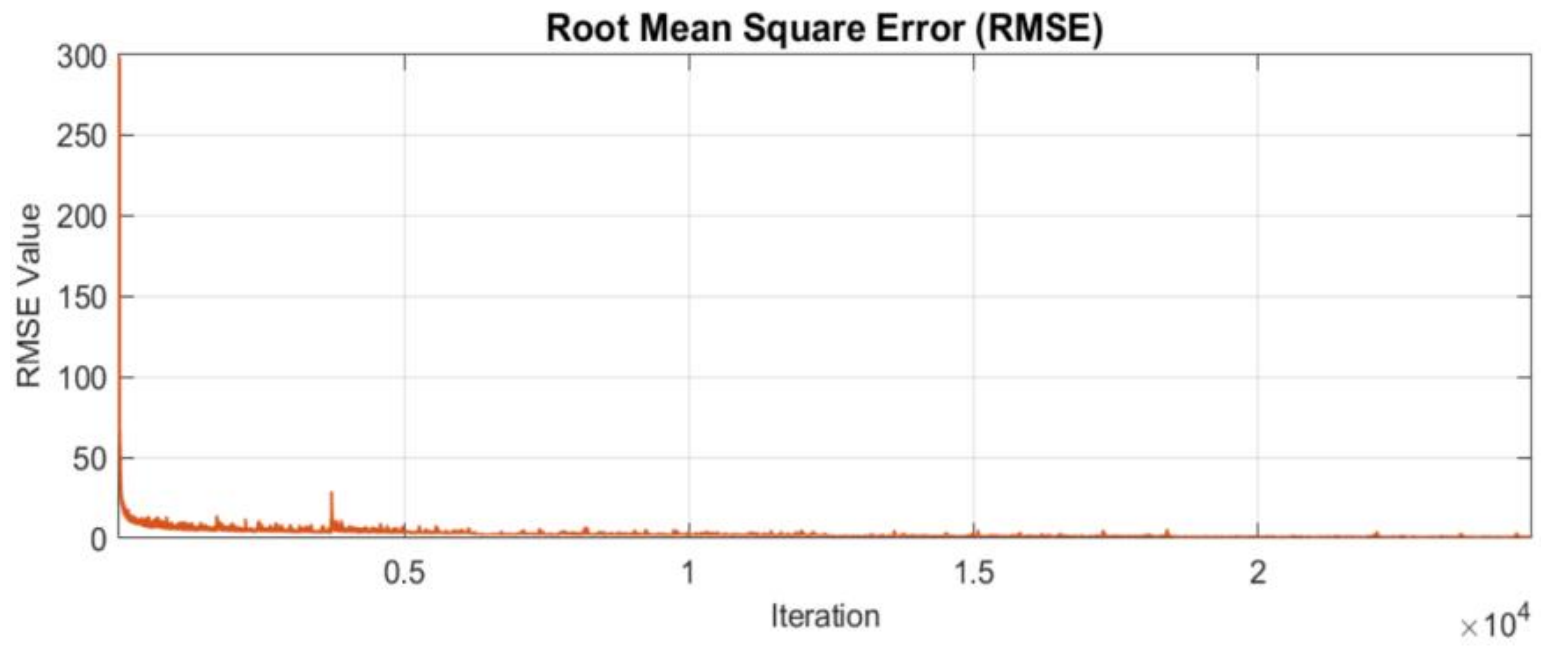

Fig. 2. Training convergence of the diffusion model (Poisson Equation)

### 4.4.3 Physics-Guided Inference and Comparative Analysis

Following the training phase, we conducted a comparative validation of the proposed method against benchmark algorithms across four test cases with distinct diffusion coefficients.

**Iterative Refinement Dynamics:** Figure 3 illustrates the full dynamics of the physics-guided denoising process as it progressively recovers the two-dimensional solution field of the Poisson equation. Driven by the physical residual guidance, the generated fields converge steadily toward the ground truth from the onset of inference. During the early phase of reverse generation (large $t$), the model leverages data priors from the network to rapidly reconstruct low-frequency contours and global topological structures. In the terminal phase (small $t$), the physics-guided mechanism fine-tunes these structures by minimizing PDE residuals. This correction effectively mitigates amplitude errors caused by distributional shifts, ensuring the final solution strictly adheres to the governing physical laws.

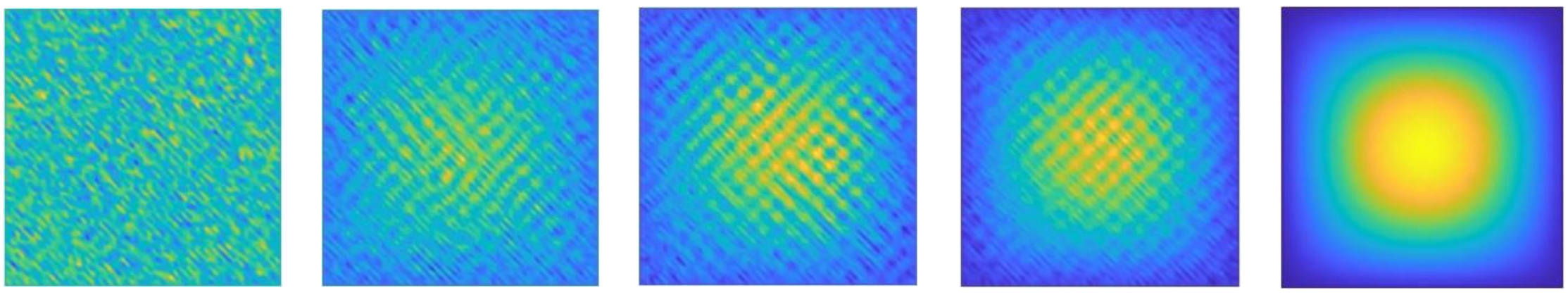

Fig. 3. Iterative evolution of the physics-guided denoising process for the Poisson equation

**Interpolation Performance:** Next, Figure 4 presents the solution results for parameters (e.g., $\kappa =$ 1.35 and 1.65) sampled within the training interval but excluded from the discrete training set, thereby evaluating the model's interpolation generalization. Quantitatively, the unguided diffusion model—operating as a purely stochastic, data-driven generator—fails to solve the Poisson equation accurately. While it generates different realizations for each run due to its inherent stochasticity, these predictions tend to gravitate toward the statistical mean of the training data, yielding high full-field relative $L_2$ errors of 49.84% and 35.49%. Error map analysis reveals that while the unguided model recovers global morphology, it fails to capture the high-frequency amplitude variations induced by parameter shifts, with predictions tending toward the statistical mean of the training dataset.

In contrast, our physics-guided framework achieves a substantial improvement in precision reliability. Despite being initialized from random noise and perturbed by stochastic fluctuations during the reverse process, the physics-guided mechanism robustly steers the denoising trajectory toward the exact solution manifold. As shown in the Fig 4, the full-field relative $L_2$ errors drop significantly to 3.43% and 4.34%, with maximum absolute errors constrained to approximately 0.008. The cross-sectional map comparison further corroborates this: the profiles from the physics-guided model nearly overlap with the ground truth, whereas the unguided model deviates significantly. These results provide strong evidence that the physics-guided mechanism effectively corrects the distributional bias of data-driven models, robustly steering the denoising trajectory toward the exact solution.

**Comparison with Physics-Informed Learning (PINNs):** Furthermore, we benchmarked our method against Physics-Informed Neural Networks (PINNs). In the experiments, the PINN was trained for 3,000 epochs, using 400 boundary points, 200 internal data points, and 1,000 collocation points. As evidenced in Fig 4, the diffusion model with physics-guided inference achieves solution accuracy comparable to, or in some instances exceeding, the PINN baseline. Crucially, while PINNs require re-sampling and complete retraining for every new parameter configuration—a computationally intensive process—our method demonstrates superior deployment efficiency. Once the diffusion prior is trained, it can generate high-precision solutions for unseen parameters in seconds via a single inference pass, eliminating the need for additional parameter updates.

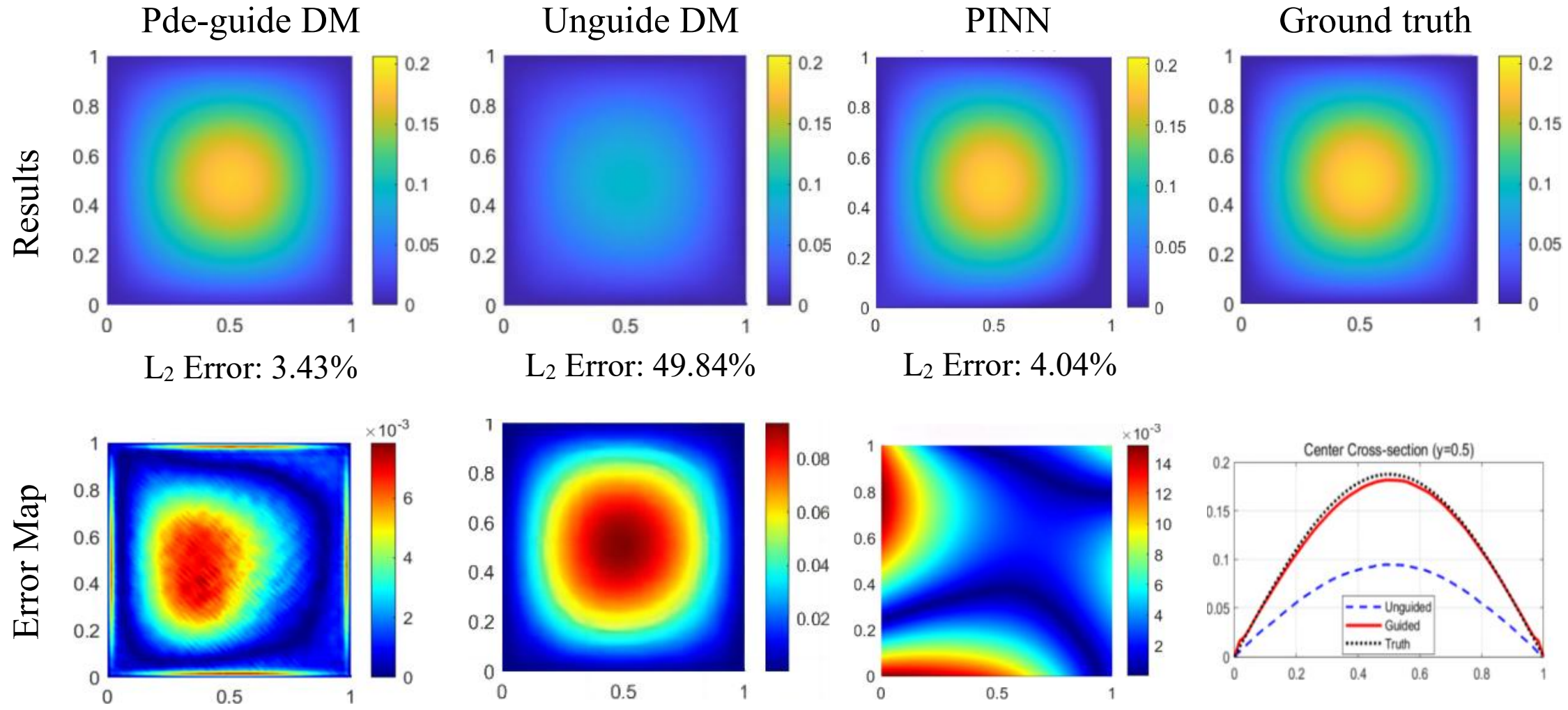

(a) Comparison of solution results with the coefficient $\kappa$=1.35

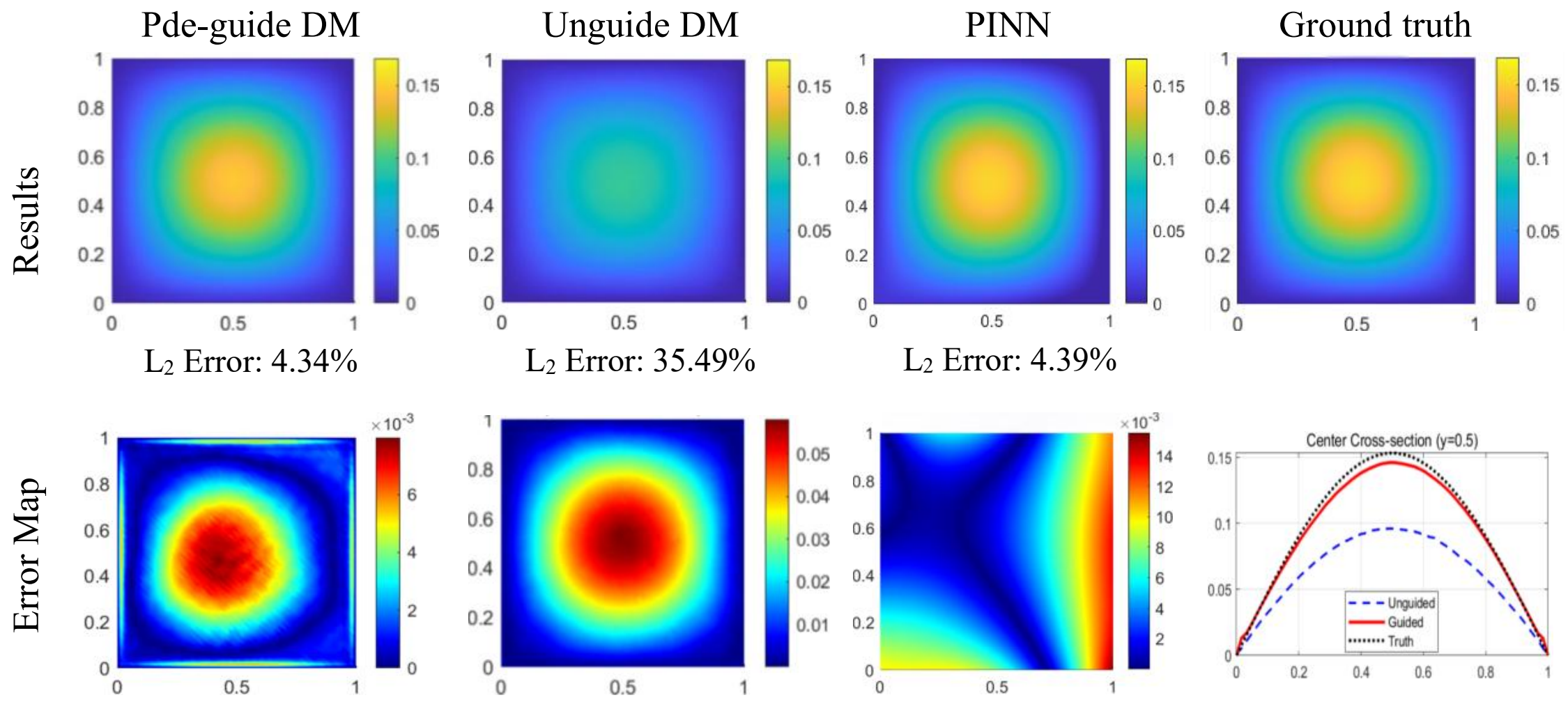

(b) Comparison of solution results with the coefficient $\kappa$=1.65

Fig. 4. Comparison of solution accuracy for the Poisson equation under interpolation coefficients. The panels display the predicted fields and corresponding error maps for the unguided and physics-guided models. The line plots illustrate the cross-sectional profiles at $y = 0.5$, demonstrating the close agreement between the physics-guided results and the ground truth.

**Extrapolation Capability:** Finally, we evaluated the model's extrapolation capability by selecting test parameters (e.g., $\kappa = 0.90$ and 2.05) that lie strictly outside the training distribution. In this regime, the unguided model exhibits severe performance degradation, as the generative prior lacks local data support at these distribution extremes. In contrast, the proposed method maintains an accuracy level consistent with its interpolation performance, robustly capturing the ground-truth amplitude variations. Notably, our framework consistently outperforms the PINN baseline (even with its case-specific retraining) in terms of solution precision under these out-of-distribution conditions. These findings demonstrate that by integrating explicit physical constraints, our method possesses reliable extrapolation capabilities, enabling robust and accurate solution generation in regions where training data is entirely absent.

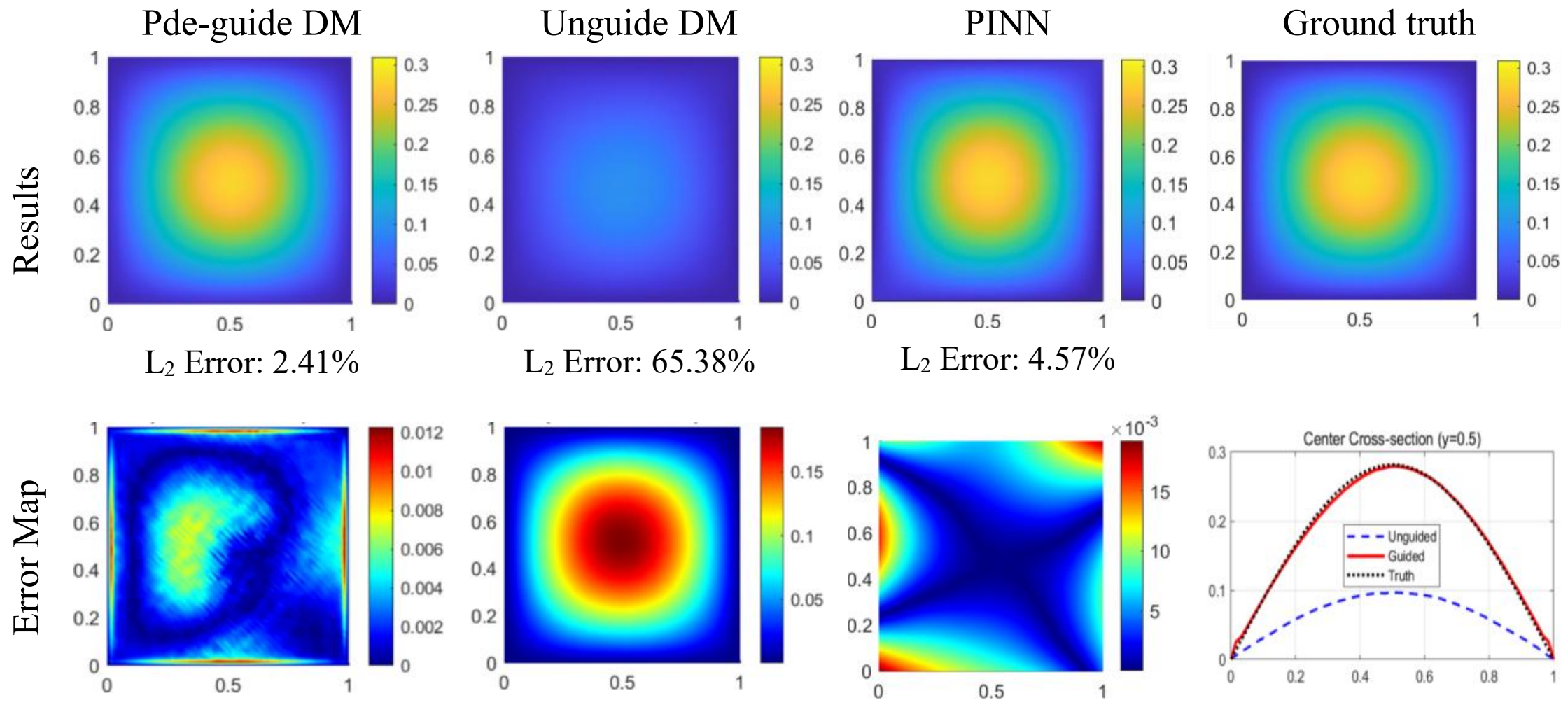

(a) Comparison of solution results with the coefficient $\kappa$=0.90

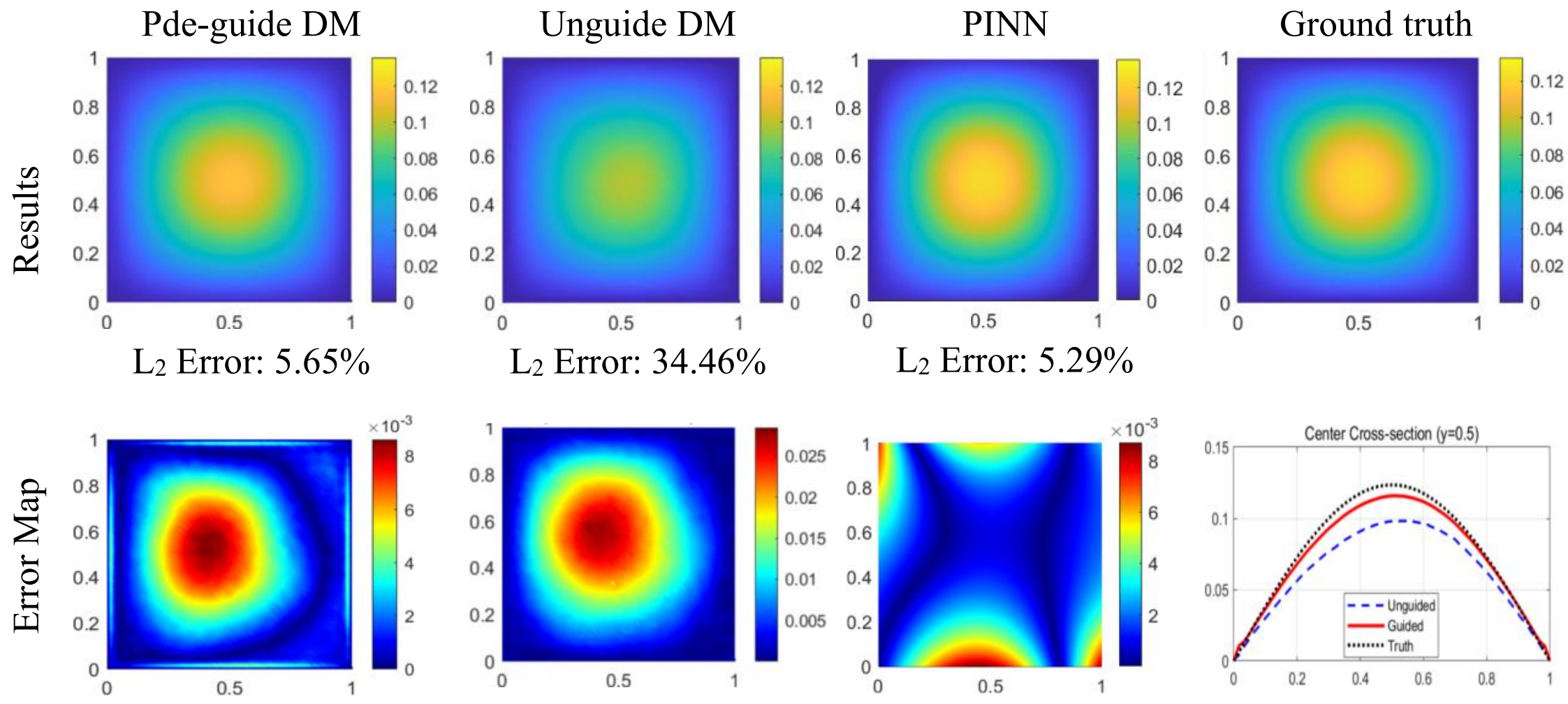

(b) Comparison of solution results with the coefficient $\kappa$=2.05

Fig. 5. Comparative analysis of solution accuracy for the Poisson equation under extrapolation coefficients. The plot configurations and performance metrics are consistent with those described in Fig. 4, demonstrating the model's robust generalization beyond the training distribution range.

## 4.5 Heat Diffusion Equation

To further validate the model's capability in solving time-dependent partial differential equations, we consider the one-dimensional heat diffusion equation.

### 4.5.1 Experimental Setup

In this experiment, the source term is set to zero to focus on the intrinsic diffusion dynamics. To evaluate the model's accuracy in capturing spatiotemporal evolution, we constructed a dataset covering a spectrum of thermal diffusion coefficients. Following our unified framework, the 1D transient problem is transformed into a 2D steady-state problem by treating time as a spatial-like dimension, thereby mapping the entire spatiotemporal evolution onto a single manifold.

The training data, consisting of 4,000 samples with coefficients uniformly sampled within the interval $[0.02, 0.05]$, were generated via a high-precision numerical solver. The network was trained for 200 epochs, resulting in 62 iterations per epoch and a total of 12,400 gradient updates based on the batch size in Table 1, employing a decay schedule that halved the learning rate every 100 epochs. For the inference phase, four distinct cases were designed to evaluate generalization: two coefficients ($\alpha = 0.031, 0.043$) within the training interval were selected to test interpolation accuracy, while two coefficients ($\alpha = 0.015, 0.055$) strictly outside the range were used to validate extrapolation capability. In these experiments, the reverse diffusion process was set to 750 steps, with the physics-guided step size calibrated to the specific residual dynamics of the heat equation.

### 4.5.2 Training Dynamics

Prior to evaluating the performance of the physics-guided inference solver, we examine the diffusion model's acquisition of spatiotemporal heat conduction patterns. By processing noisy spatiotemporal fields, the model learns to reconstruct characteristic thermal diffusion features across varying noise intensities. Figure 6 presents the training convergence, where the red line tracks the Root Mean Squared Error (RMSE). The loss functions exhibit a rapid initial decline before plateauing after approximately 6,000 iterations. This suggests that the U-Net has effectively captured the intrinsic governing laws of the heat equation within the transformed 2D spatiotemporal domain, establishing a robust prior for the subsequent steady-state inference of transient evolution.

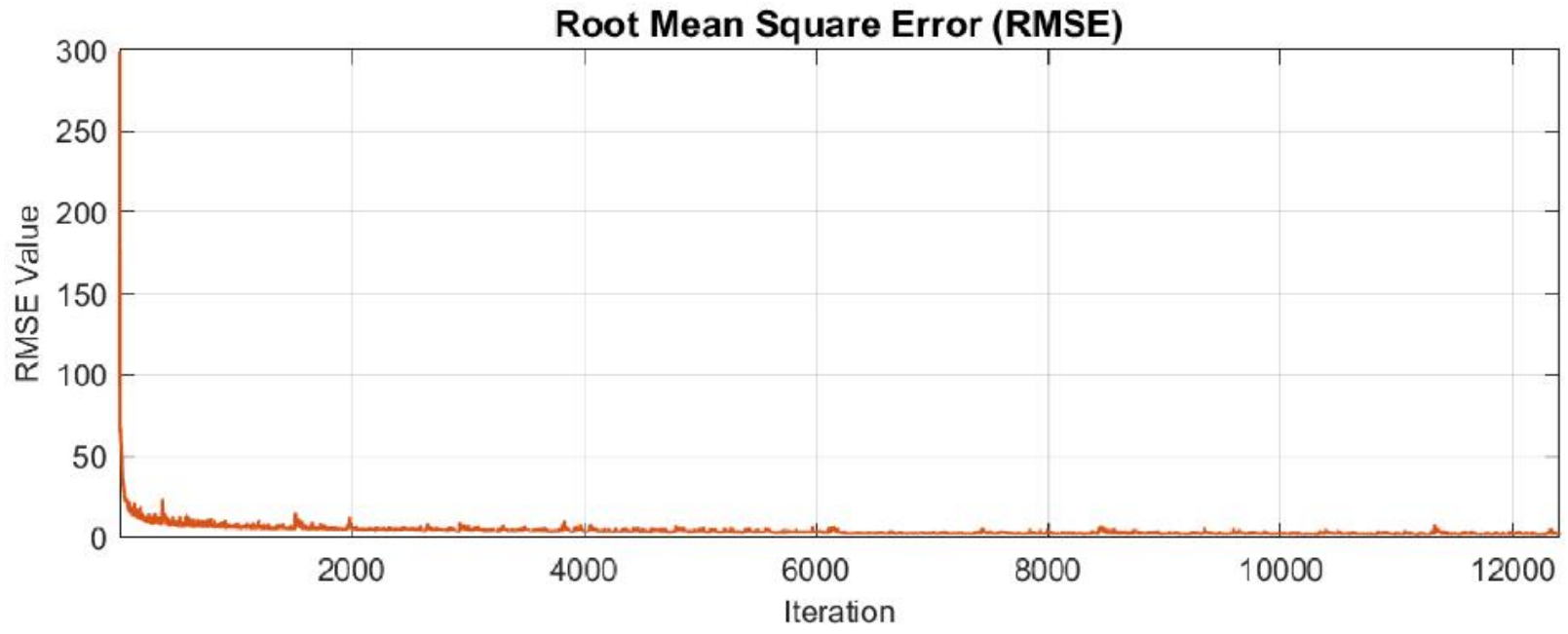

Fig. 6. Training convergence of the diffusion model for the spatiotemporal heat equation

### 4.5.3 Physics-Guided Inference and Comparative Analysis

Following the training phase, we evaluated the proposed method against benchmark algorithms using test cases with varied thermal diffusion coefficients.

**Iterative Refinement Dynamics:** Figure 7 illustrates the dynamic denoising process of the diffusion model with physics-guided inference for the time-dependent heat equation. Starting from pure Gaussian noise, the model progressively reconstructs the spatiotemporal temperature field. During the early stage of the reverse generation (large $t$), the data prior, augmented by physical guidance, rapidly establishes the global trend of thermal decay. In the terminal stage ($t \to 0$), the denoising process suppresses irregular oscillations while the physics-guided mechanism fine-tunes the amplitude and energy diffusion profiles. This ensures that the generated field captures the characteristic transition from the initial transient state to a stabilized distribution.

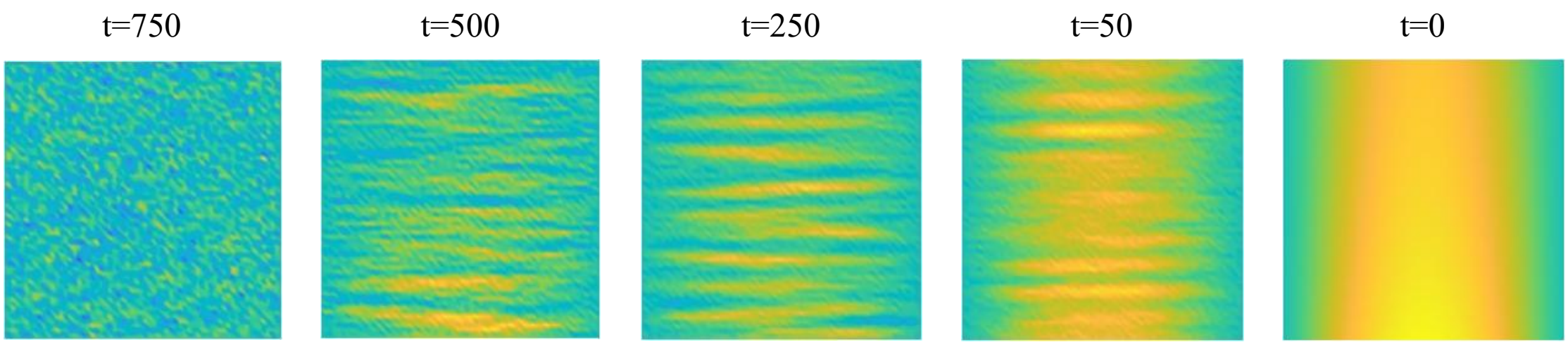

Fig. 7. Iterative evolution of the physics-guided denoising process for the spatiotemporal heat equation

**Interpolation Performance:** Figure 8 evaluates the model's interpolation capability by presenting solution fields for coefficients that, while remaining within the training bounds, represent unseen coordinates in the parameter space. This test verifies the framework's precision in interpolating the underlying spatiotemporal dynamics between discrete training points. Quantitative analysis shows that the purely data-driven unguided model struggles to differentiate between varying diffusion rates; its stochastic sampling lacks physical grounding, causing predictions to regress toward the statistical mean of the training data. This leads to high full-field $L_2$ errors of 20.27% and 18.22%, respectively. In contrast, our proposed method significantly mitigates these discrepancies by incorporating physical energy constraints. Notably, despite starting from random initial noise, the physics-guided mechanism ensures that the denoising process converges accurately and stably to the target solution. Although minor deviations persist at local extrema—with a maximum absolute error of approximately 0.06, the overall full-field L2 error remains consistently between 3% and 4%. This demonstrates the model's ability to accurately and robustly resolve the dynamic heat conduction process on unseen parameters within the training manifold.

**Comparison with Physics-Informed Learning (PINNs):** We Furthermore compared the proposed method with the PINN baseline. In the experiment, the PINN was trained for 3,000 epochs, using the same amount of training data as that of in the Poisson equation experiment. The results in Figure 8 show that the PINN achieves a stable full-field $L_2$ error below 2%, our generalizable solver maintains a competitive performance gap of only approximately 2%. However, it is critical to highlight the trade-off in computational efficiency: PINNs require re-sampling and a complete re-optimization cycle for every new diffusion coefficient. In contrast, our method performs zero-shot inference; once trained, it provides a reliable full-field solution for unseen coefficients within seconds, bypassing the need for computationally expensive retraining.

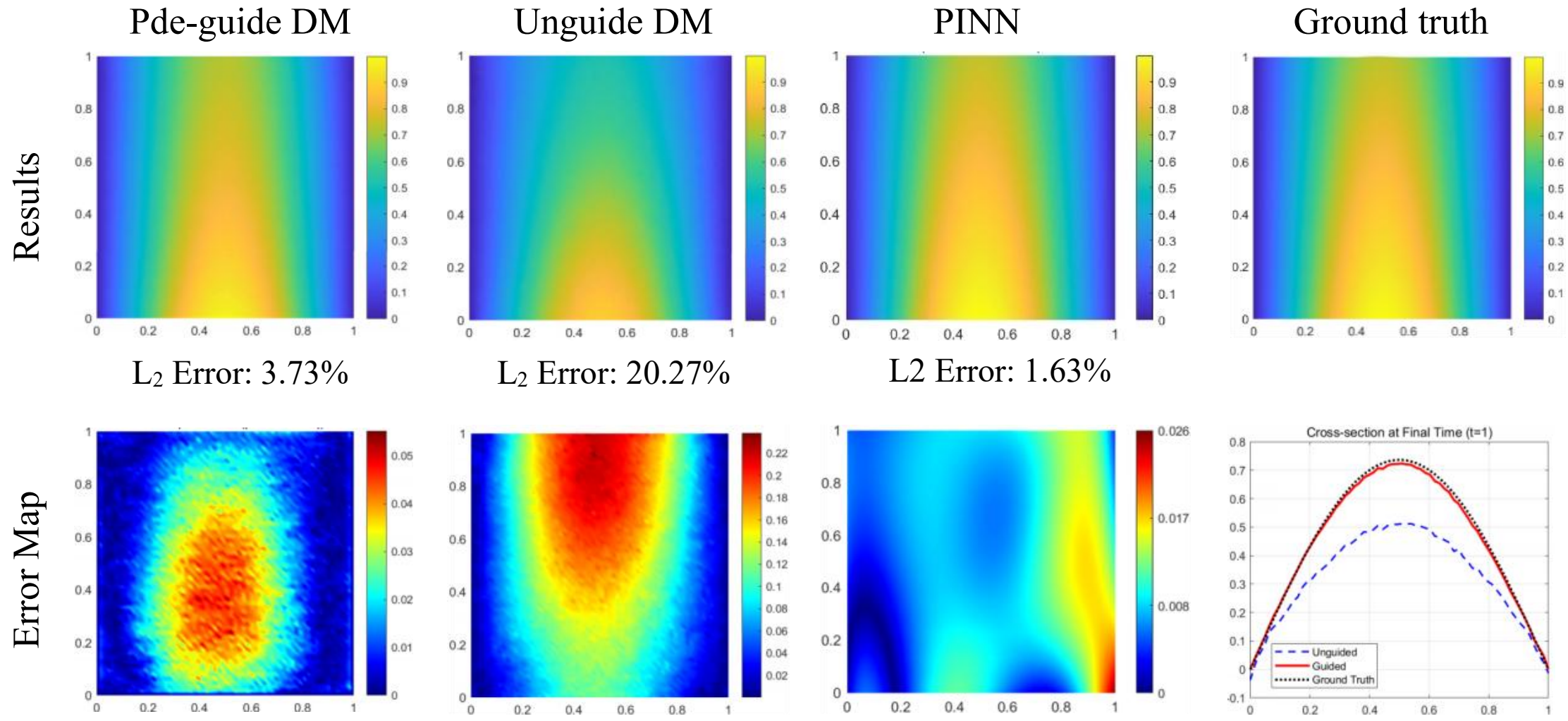

(a) Comparison of solution results with the thermal diffusion coefficient α=0.031

Pde-guide DM Unguide DM PINN Ground truth

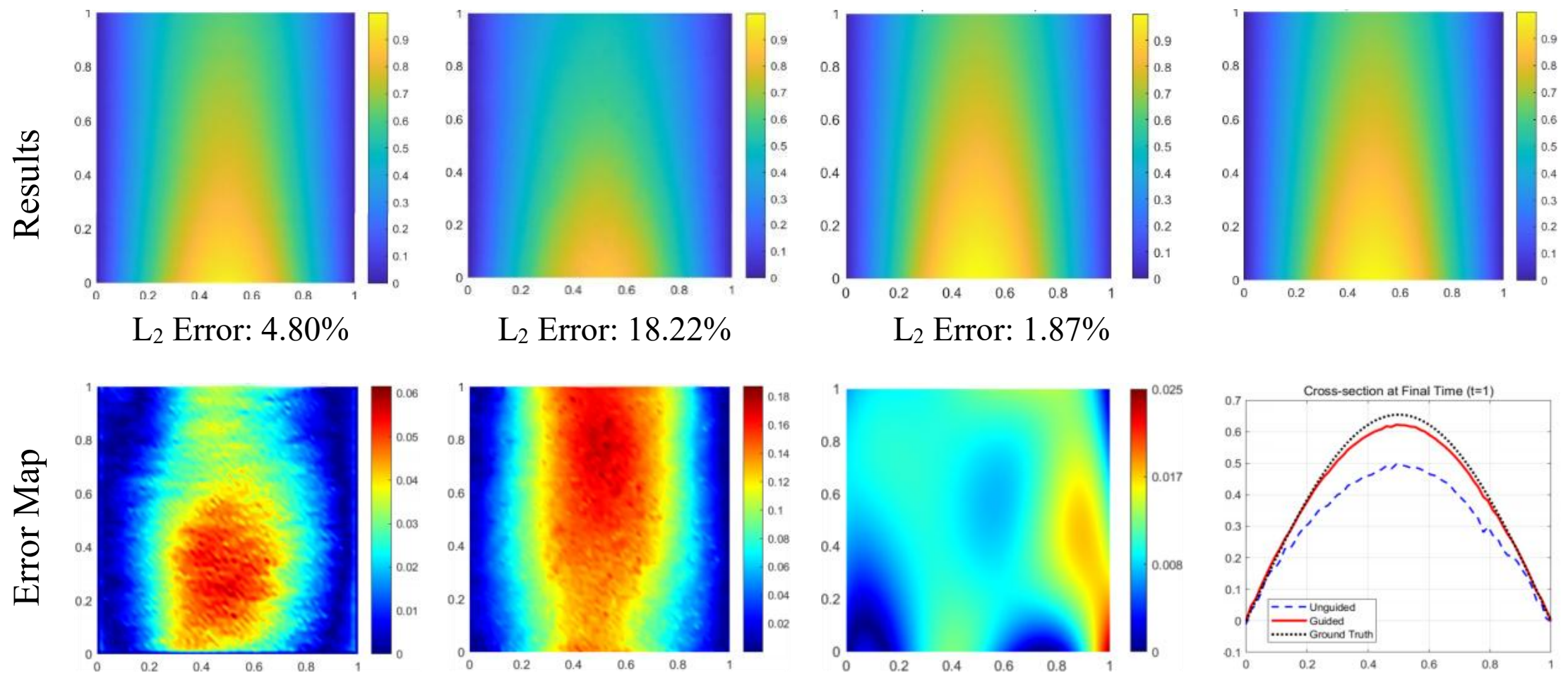

(b) Comparison of solution results with the thermal diffusion coefficient α=0.043

Fig. 8. Comparison of solution accuracy for the heat equation under interpolation coefficients. The panels display the spatiotemporal evolution from the initial transient state to steady-state decay. The plot configurations are consistent with Fig. 4, with line plots representing the temporal slice at $t = 1.0$.

**Extrapolation Capability:** Finally, we assessed the model's robustness beyond the training distribution using coefficients strictly outside the original interval. While the unguided model fails to produce physically meaningful solutions, the proposed method maintains an accuracy level consistent with its interpolation results. Despite the slight precision gap compared to a specifically retrained PINN, our model successfully generates temperature fields that strictly obey the underlying physical laws. This provides strong evidence that the physics-guided mechanism enables robust extrapolation for time-dependent parabolic PDEs, even in the absence of local training data.

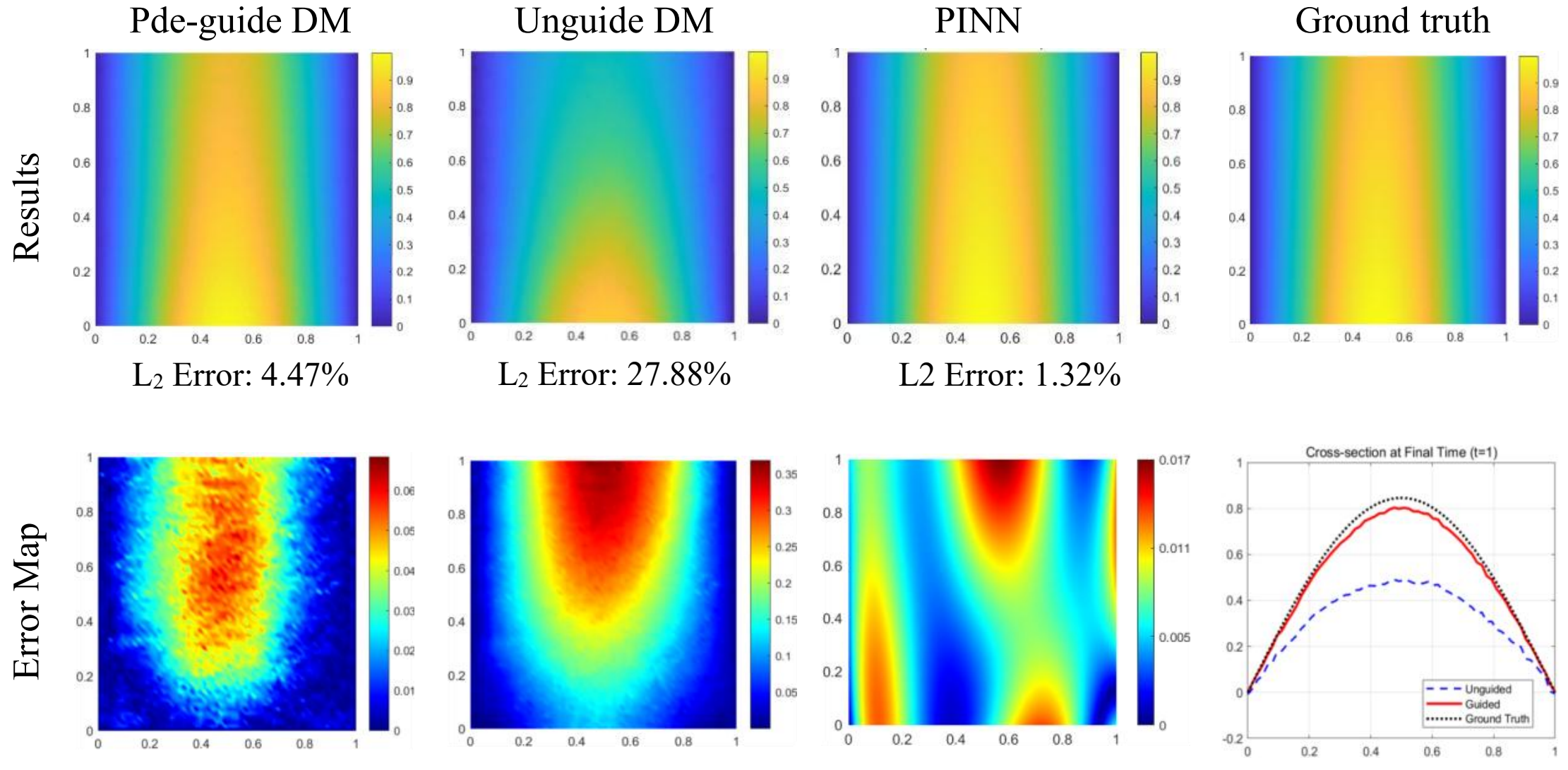

(a) Comparison of solution results with the thermal diffusion coefficient α=0.015

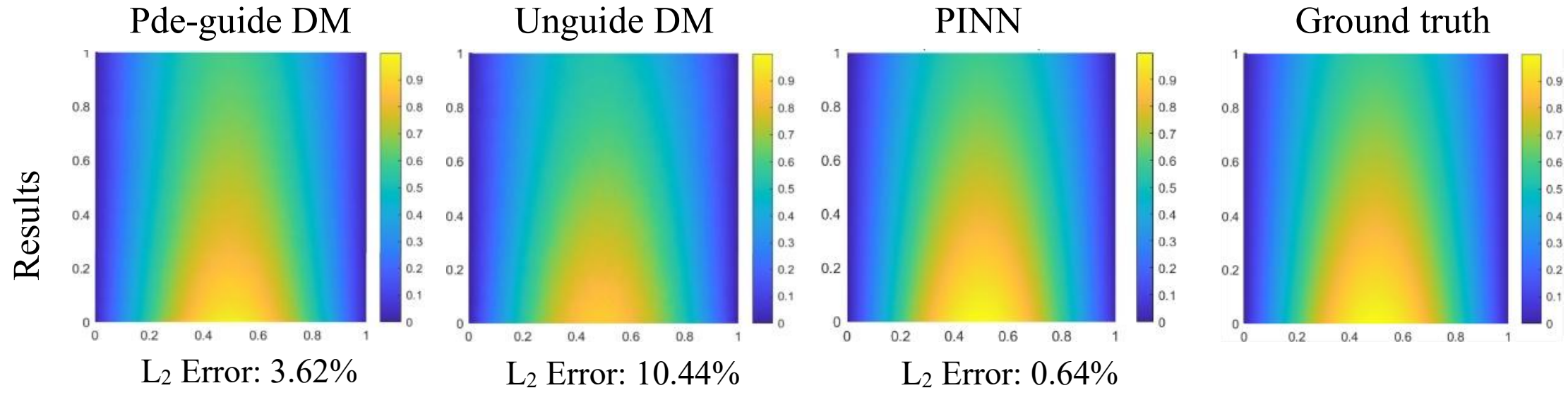

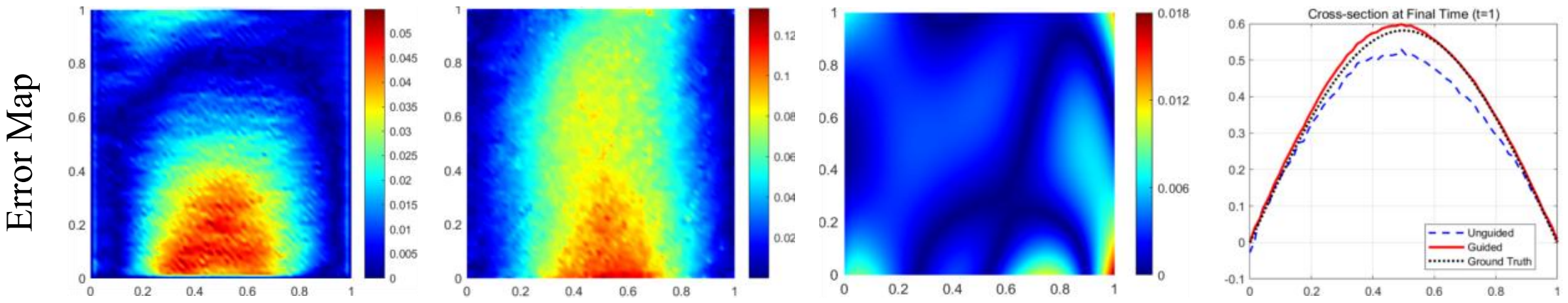

(b) Comparison of solution results with the thermal diffusion coefficient α=0.055

Fig. 9. Comparative analysis of solution accuracy for the heat equation under extrapolation coefficients. The plot configurations and performance metrics are consistent with those described in Fig. 8, illustrating the model's robust tracking of the transient-to-steady thermal evolution even beyond the training distribution range.

## 4.6 Burgers Equation

To examine the model's performance on nonlinear advection-diffusion systems and its capability to handle shock wave and resolve sharp gradients, we consider the one-dimensional Burgers' equation.

### 4.6.1 Experimental Setup

To rigorously evaluate the framework's accuracy in capturing nonlinear spatiotemporal evolution and shock formation, we constructed a dataset with varied viscosity coefficients. Similar to the heat equation setup, the 1D transient velocity field is mapped onto a 2D spatiotemporal manifold. The training data, comprising 4,000 samples with coefficients uniformly sampled within the interval $[0.01, 0.03]$, were generated via a high-precision numerical solver. The model was trained for 300 epochs, resulting in 62 iterations per epoch and a total of 18,600 gradient updates based on the batch size in Table 1, employing a decay schedule that halved the learning rate every 100 epochs. For the inference phase, four distinct cases were designed: two coefficients ($\nu = 0.017, 0.024$) within the training interval were selected to assess interpolation generalization, while two others ($\nu = 0.0075, 0.0325$) strictly outside the range were used to validate extrapolation robustness. The reverse diffusion process was set to 750 steps, with the physics-guided step size tailored to the nonlinear residual characteristics of the Burgers' equation.

### 4.6.2 Training Dynamics

Prior to evaluating the physics-guided solver, we analyze the model's acquisition of the nonlinear spatiotemporal patterns inherent in the Burgers' equation. By processing noisy fields across various noise intensities, the U-Net learns to reconstruct both the sharp shock fronts and the smooth regions of the flow field. Figure 10 illustrates the training convergence, where the red line tracks the RMSE. The loss exhibits a rapid initial decline and reaches a stable plateau after approximately 8,000 iterations. This convergence indicates that the network has effectively captured the intrinsic advection-diffusion mechanisms, establishing a high-fidelity generative prior for resolving nonlinear evolution.

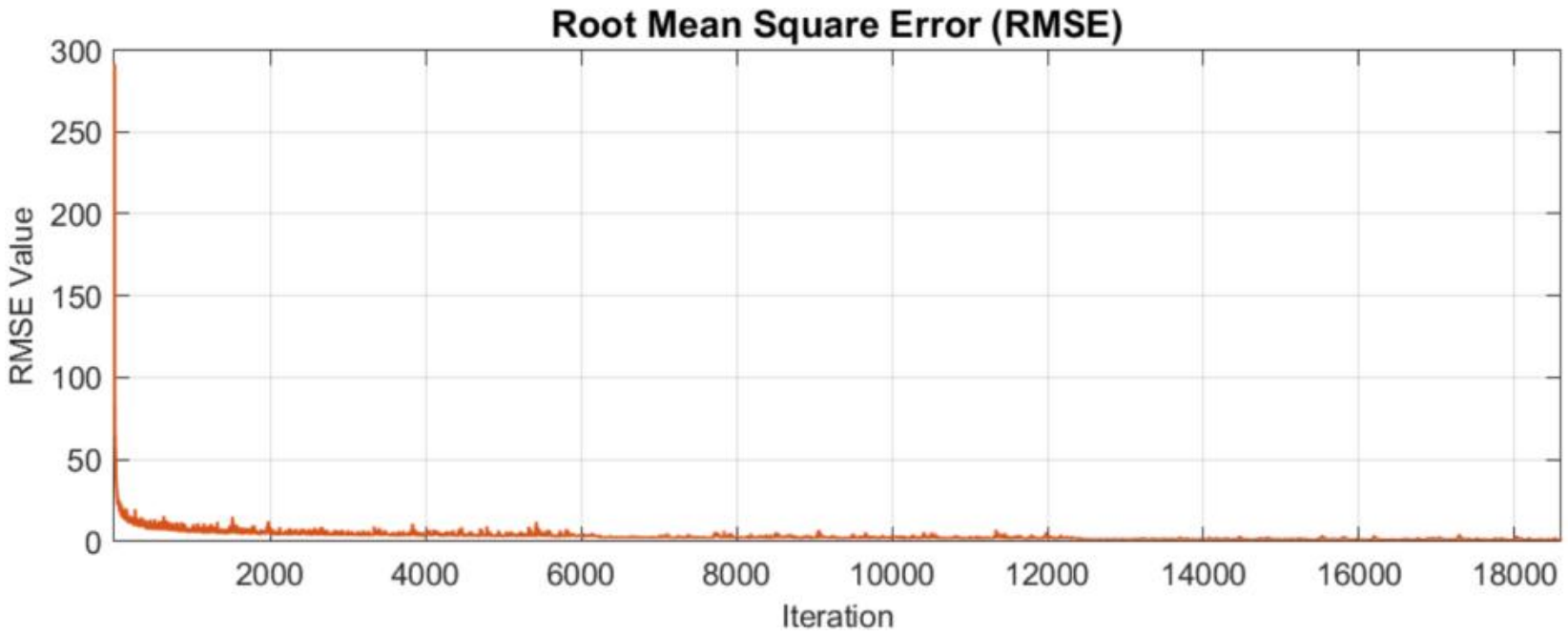

Fig. 10. Training convergence of the diffusion model for the nonlinear Burger’s equation

### 4.6.3 Physics-Guided Inference and Comparative Analysis

Following the heat equation validation, we further evaluate the model's capacity to resolve nonlinear advection-diffusion systems and capture shock wave dynamics. We selected varied viscosity coefficients to benchmark the proposed method against baseline solvers.

**Iterative Refinement Dynamics:** Figure 11 illustrates the dynamic inference process of the physics-guided diffusion model for solving the Burgers equation. Starting from pure Gaussian noise, the model progressively reconstructs the spatiotemporal evolution of shock formation and propagation. During the early reverse generation stage (large $t$)，the data prior and physical guidance collaboratively establish the global convection trends and the approximate location of the shock front. In the terminal stage ($t \rightarrow 0$), the denoising process suppresses stochastic oscillations while the physics-guided mechanism refines the steep gradients near the shock and corrects the velocity amplitudes, ensuring high-fidelity recovery of nonlinear features..

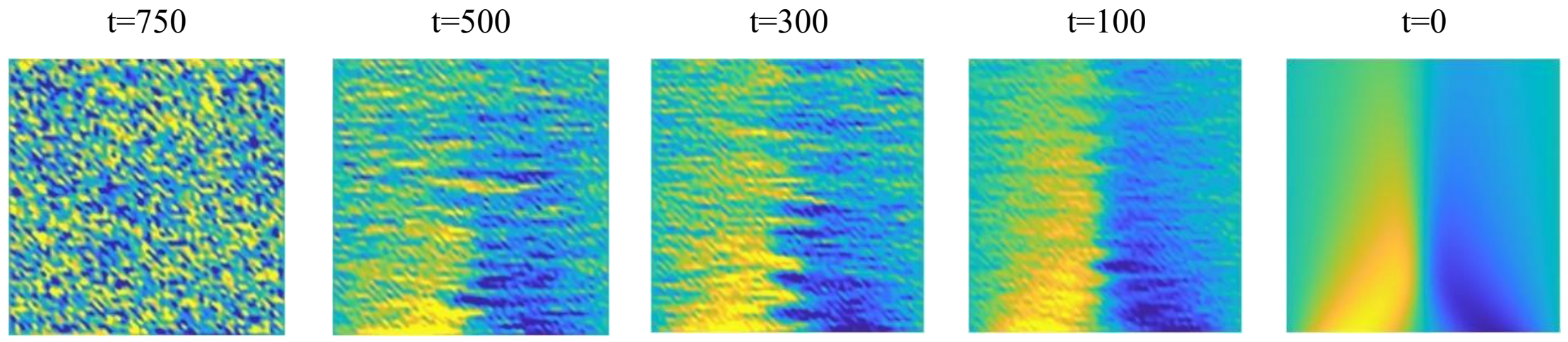

Fig. 11. Iterative evolution of the physics-guided denoising process for the nonlinear Burgers' equation

**Interpolation Performance:** Figure 12 presents the results for the Burgers equation on the interpolation test set. Quantitative analysis reveals that the unguided model fails to accurately capture the steepness and position of shock waves across different viscosity coefficients, as its predictions tend to regress toward the training ensemble's mean. In contrast, our method effectively mitigates these errors by incorporating the nonlinear physical residuals of the Burgers' equation. Crucially, while the reverse process initiates from random Gaussian noise, the physics-guided mechanism consistently anchors the stochastic trajectory, ensuring stable convergence to the physically accurate solution. Although the high-frequency gradients at the shock front pose a significant challenge, resulting in minor persistent deviations at the extrema, the overall full-field L2 errors are constrained to 1.89% and 3.39%, with absolute error below 0.04. This proves that the framework can accurately reconstruct complex nonlinear dynamics for unseen parameters within the training manifold.

**Comparison with Physics-Informed Learning (PINNs):** We further compared our approach with the PINN baseline. In this experiment, the PINN was trained for 3,000 epochs, with other settings consistent with previous experiments. As shown in Figure 12, while the PINN achieves a superior accuracy (approximately 1% $L_2$ error), our method functions as a generalizable solver. Despite a marginal precision gap (1% – 2%), our framework bypasses the need for the computationally expensive retraining required by PINNs for every new coefficient. For unseen scenarios, it delivers a physically consistent full-field solution in seconds, achieving an optimal balance between generalization and efficiency.

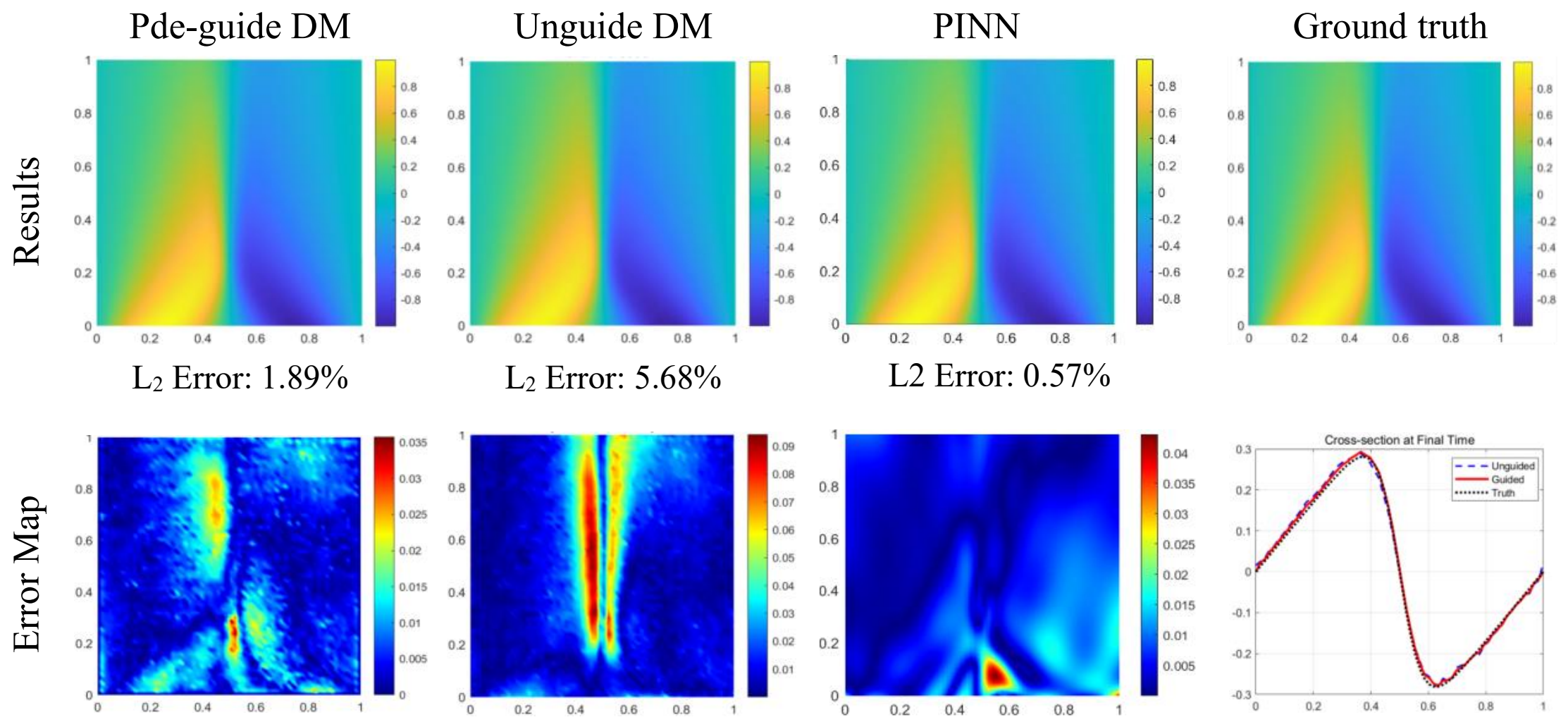

(a) Comparison of solution results with the kinematic viscosity ν=0.017

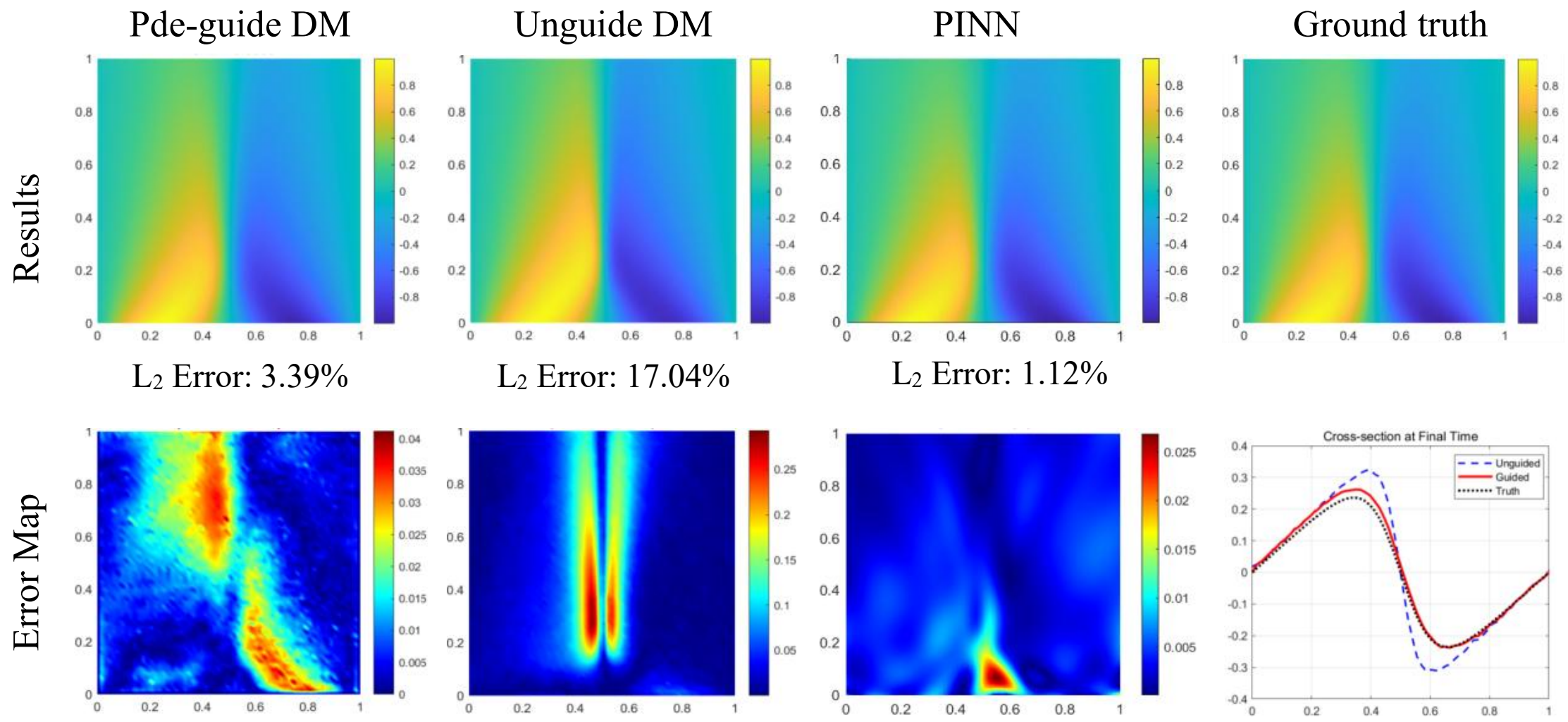

(b) Comparison of solution results with the kinematic viscosity ν=0.024

Fig. 12. Comparison of solution accuracy for the Burgers' equation under interpolation coefficients. The panels display the reconstruction of nonlinear shock wave propagation. The plot configurations are consistent with Fig. 4, with line plots representing the temporal slice at the final state $t = 1.0$.

**Extrapolation Capability:** Finally, we assessed the model's robustness under extreme, out-of-distribution conditions, specifically focusing on cases with lower viscosity coefficient $\nu = 0.0075$ where convection dominates. In such scenarios, the shock front becomes significantly steeper than any sample in the training set. The unguided model, limited by the training distribution, produces an overly smoothed solution that fails to resolve the sharp gradient.

As illustrated in Fig 13(a), the physics-guided mechanism forces the model to accurately resolve the discontinuity and steep gradient features at the shock location. Even when corresponding training data is absent, the model generates flow fields that strictly adhere to the underlying physical laws. Figure 13(b) further confirms the model's stability in solving cases with relatively smooth fronts. These results demonstrate that the physics-guided diffusion model possesses robust extrapolation capabilities for nonlinear PDEs characterized by shock waves.

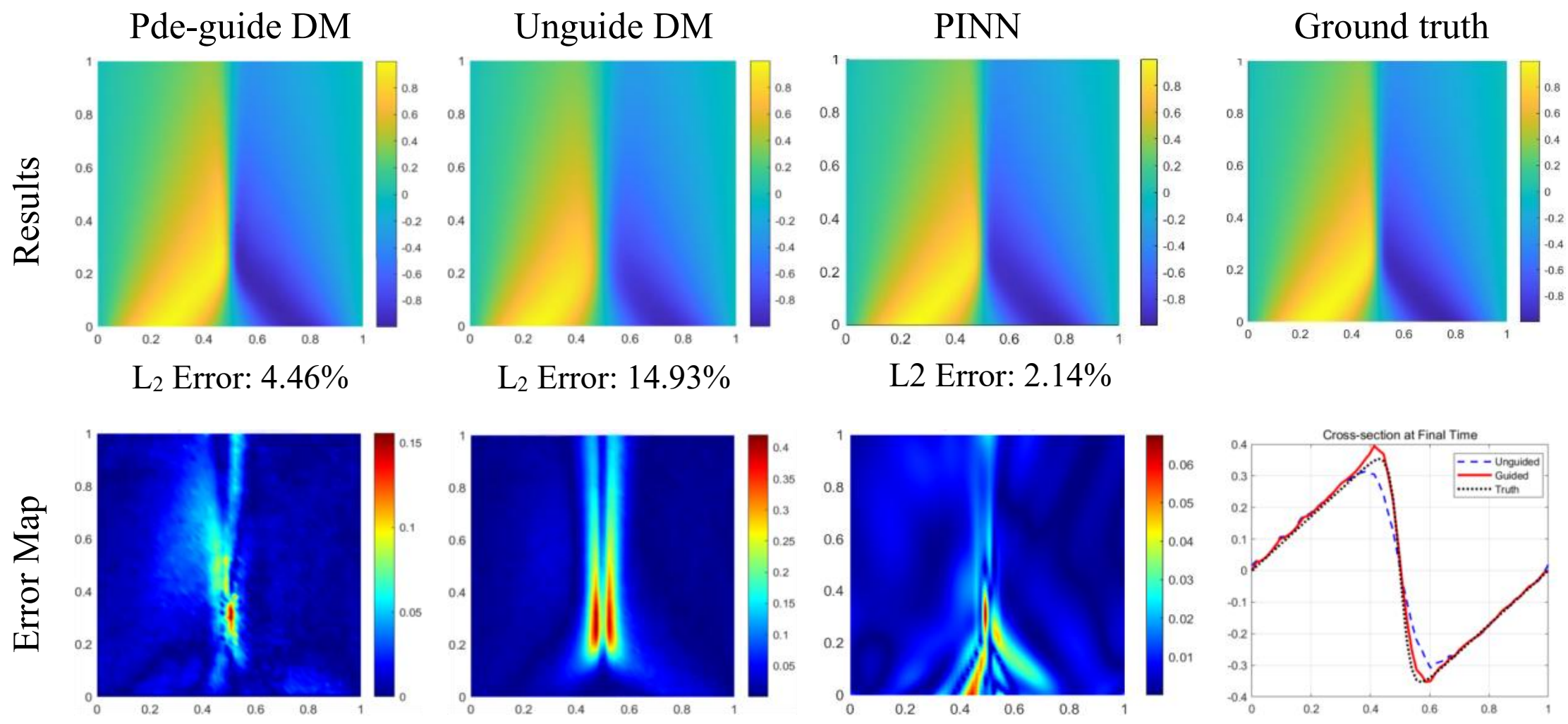

(a) Comparison of the solution with the kinematic viscosity ν=0.0075

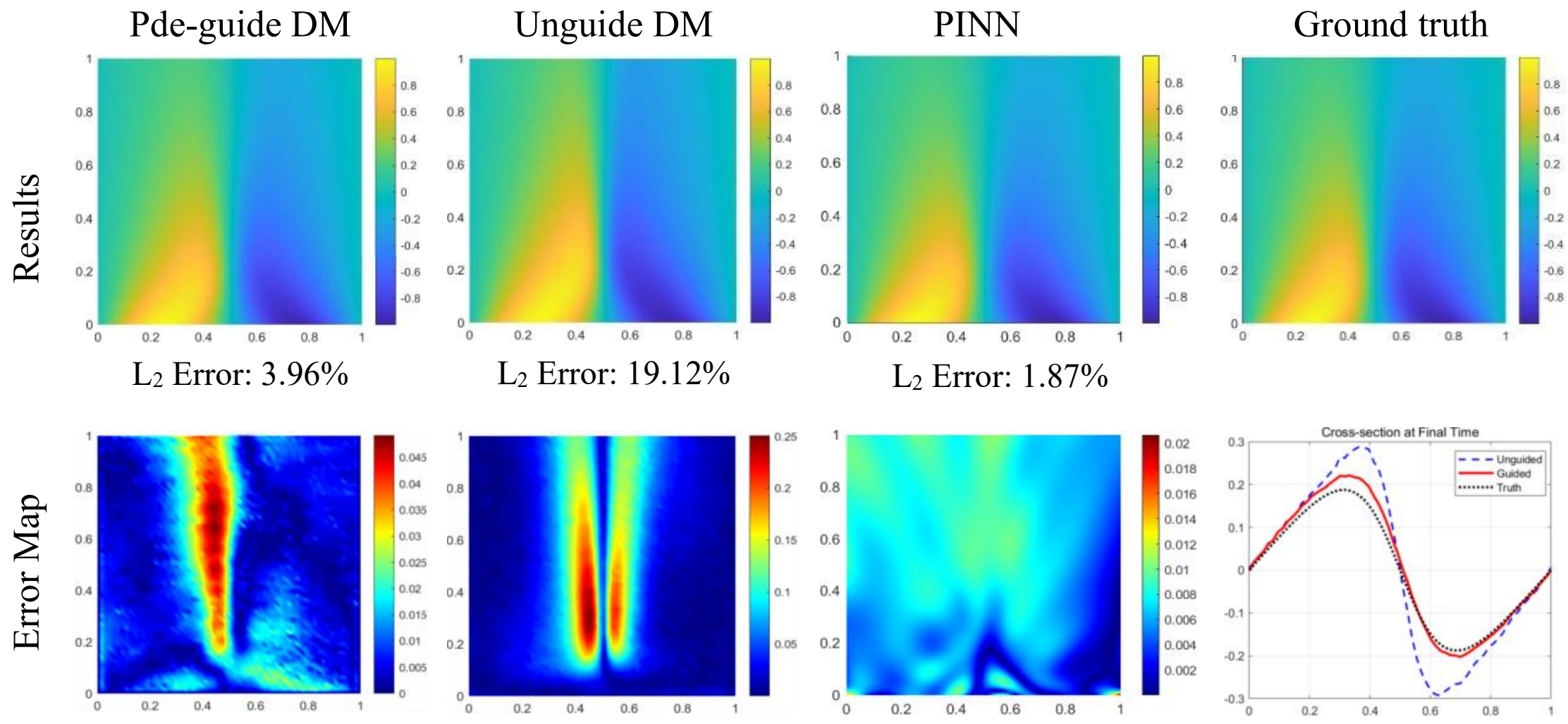

(b) Comparison of the solution with the kinematic viscosity ν=0.0325

Fig. 13. Comparative analysis of solution accuracy for the Burgers' equation under extrapolation coefficients. The plot configurations and metrics are consistent with Fig.12, illustrating the model's robust capability to resolve extremely steep gradients and nonlinear discontinuities beyond the training distribution range.

## 4.7 Quantitative Analysis and Discussion

### 4.7.1 Statistical Accuracy and Robustness

To evaluate the robustness of the proposed model, we conducted fifty independent inference runs for four representative coefficients across the three governing equations. The mean and standard deviation of the full-field relative $L_2$ errors were calculated and are summarized in Table 2. Under consistent experimental configurations, the standard deviation of the $L_2$ error remained consistently low, ranging from 0.5% to 1.5% across various physical coefficients. This suggests that the physics-guided energy constraints effectively anchor the diffusion inference process onto a unique physical manifold, thereby suppressing the stochastic drift inherent in vanilla generative models. Combined with the low average error, these findings confirm that our model is not only robust against random sampling interference but also maintains a high degree of fidelity in both interpolation and extrapolation tasks, satisfying the rigorous requirements of scientific computing.

Table 2. Statistical robustness of the physics-guided diffusion model across multiple inference trials (Mean and standard deviation are in percentage).

| Types \ Method | Poisson Equation | Heat Diffusion Equation | Burger Equation |
|---|---|---|---|
| Coefficient 1 | $\kappa = 1.35$<br>$2.8425 \pm 0.8431$ | $\alpha = 0.031$<br>$4.3500 \pm 0.8374$ | $\nu = 0.017$<br>$3.1160 \pm 1.3245$ |
| Coefficient 2 | $\kappa = 1.65$<br>$3.7160 \pm 0.9401$ | $\alpha = 0.043$<br>$3.8725 \pm 1.4023$ | $\nu = 0.024$<br>$2.9625 \pm 0.9789$ |
| Coefficient 3 | $\kappa = 0.90$<br>$2.4400 \pm 0.1402$ | $\alpha = 0.015$<br>$4.8640 \pm 1.4799$ | $\nu = 0.0075$<br>$4.6025 \pm 0.7649$ |
| Coefficient 4 | $\kappa = 2.05$<br>$4.7398 \pm 1.0561$ | $\alpha = 0.055$<br>$4.0300 \pm 1.2232$ | $\nu = 0.0325$<br>$5.0425 \pm 1.0821$ |

### 4.7.2 Ablation Study: Impact of Gaussian Smoothing

To verify the necessity of Gaussian smoothing within the physics-guided inference framework, we conducted an ablation study using the Poisson equation ($\kappa = 0.90$). The intermediate states of the diffusion process are characterized by significant high-frequency noise. Directly applying physics guidance involving differential operators (such as the Laplacian) to such noisy fields often triggers numerical gradient explosion, rendering the guidance mechanism ineffective. Consequently, we compared the inference performance "without Gaussian smoothing" against that "with Gaussian smoothing" ($\sigma = 0.90$) configuration under identical guidance conditions.

By progressively increasing the guidance step size, we identified the maximum physical guide time step ($\Delta t$) that maintains a stable solution to the Poisson equation: The model without smoothing was restricted to $\Delta t = 7.0\times10^{-5}$, whereas the inclusion of Gaussian smoothing allowed for a much larger stable step size of $\Delta t = 6.7\times10^{-4}$. The corresponding results are illustrated in Fig 14. As shown, the model utilizing Gaussian smoothing achieves accurate solutions at the optimized step size. In contrast, the un-smoothed model suffers from numerical instability at the same step size ($\Delta t = 6.7\times10^{-4}$); even at its own maximum stable limit ($\Delta t = 7.0\times10^{-5}$), its accuracy remains inferior to the smoothed counterpart.

However, it is important to note that while Gaussian smoothing is indispensable for suppressing noise-induced instability, excessive filtering may inadvertently blur sharp physical features, thereby limiting the final solution accuracy. To further balance numerical stability with predictive precision, adaptive filtering strategies can be implemented to handle this trade-off. Such mechanisms would dynamically adjust the smoothing intensity based on the noise level at each diffusion step, ensuring robust convergence without sacrificing the high-fidelity details of the physical field.

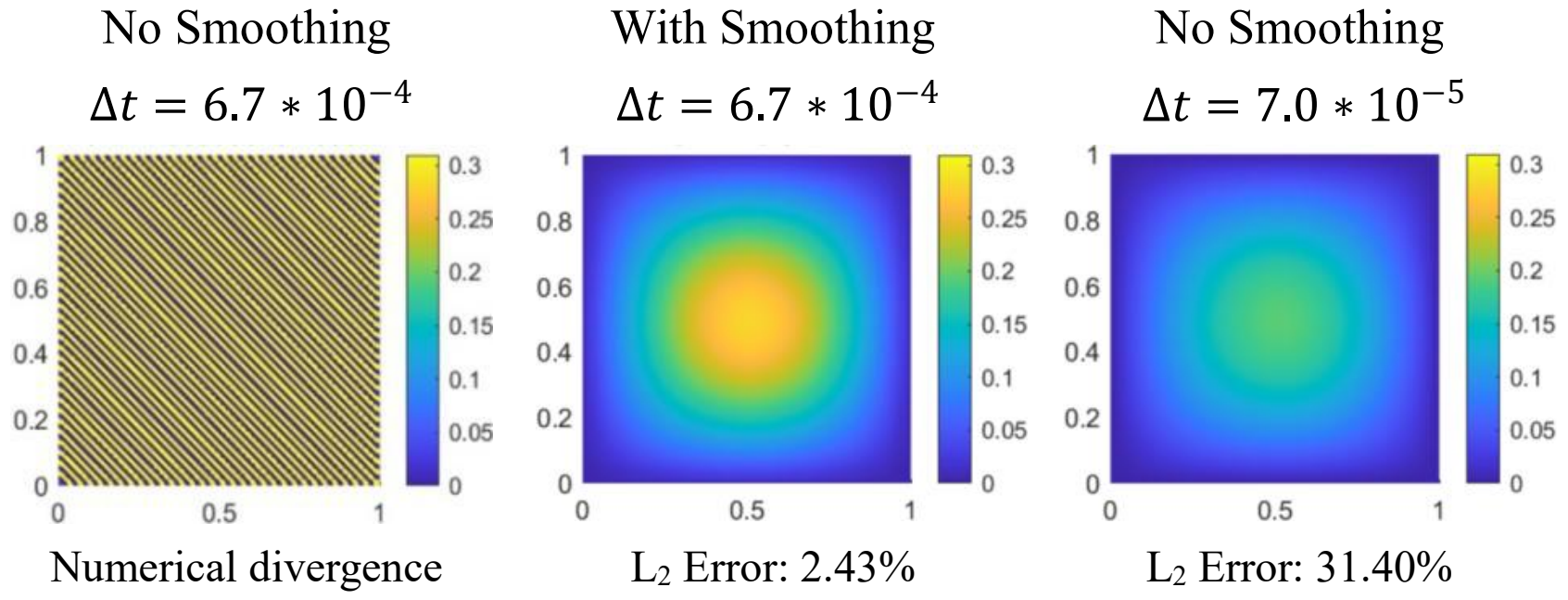

Fig. 14. Impact of Gaussian smoothing on the stability of physics-guided inference

### 4.7.3 Efficiency and General Discussion

All experiments were conducted on an NVIDIA RTX 4060 GPU. Regarding computational cost, the PINN baseline requires approximately 5 minutes of retraining to optimize parameters for each new physical coefficient. While the PINN achieves high accuracy, its total computational overhead scales linearly with the number of test cases. In contrast, our method requires only a one-time training session of 30 to 60 minutes. Once trained, the model performs zero-shot inference, generating solutions for unseen coefficients in 5 to 10 seconds without any weight updates. This significant reduction in marginal cost makes our method highly advantageous for real-time applications and scenarios requiring repeated queries across a continuous parameter space.

Comprehensive results across the Poisson, heat, and Burgers' equation experiments demonstrate that the proposed method achieves strong generalization without retraining. While purely data-driven models often regress toward the ensemble average when applied to unseen coefficients, physics-guided inference enforces the governing equations as hard constraints, effectively correcting errors caused by distribution shifts. Furthermore, despite the stochastic nature of the reverse generation process, we observed high stability across all test cases. This confirms that physics guidance provides robust gradient directions, ensuring consistent convergence to the correct physical manifold regardless of equation types or boundary conditions. Compared to baselines like PINNs and traditional numerical solvers, our method delivers comparable accuracy with superior deployment efficiency, satisfying the demands of practical scientific computing.

# 5. Conclusions

This work presents a diffusion model with physics-guided inference for solving partial differential equations, providing a principled alternative to both classical numerical solvers and physics-informed learning approaches. In contrast to existing physics-informed diffusion or PINN-based methods that embed governing equations into the training process, the proposed framework decouples learning from physics: the diffusion model is trained purely in a data-driven manner, while physical laws are enforced exclusively during the reverse inference stage through a PDE energy function, Gaussian regularization, and explicit boundary enforcement.

From a theoretical standpoint, the proposed inference dynamics can be interpreted as a stochastic gradient flow driven by the PDE energy. The resulting process converges in distribution to the Gibbs measure induced by the governing equation, and rigorously recovers classical deterministic solvers in the zero-noise limit. This establishes a clear and unified connection between diffusion models, energy-based formulations, stochastic differential equations, and traditional PDE solvers, bridging generative modeling and numerical analysis within a single mathematical framework.

From a numerical perspective, the method exhibits robust convergence from random initializations, stable behavior under stochastic perturbations, and high accuracy across elliptic, parabolic, and nonlinear PDEs. Unlike physics-informed neural networks, the proposed approach demonstrates strong generalization to unseen coefficients without retraining, owing to the explicit enforcement of physical constraints during inference rather than during learning. These results highlight the effectiveness of physics-guided diffusion as a universal inference mechanism for PDE solutions.

The proposed framework opens several promising directions for future research. First, extending the method to multi-physics systems involving coupled PDEs. Second, incorporating complex geometries, irregular domains, and heterogeneous boundary conditions will further enhance applicability to realistic engineering scenarios. Third, scaling the approach to three-dimensional large-scale problems and high-resolution industrial simulations remains an important challenge. Finally, further theoretical investigation into convergence rates, error estimates, and noise–resolution trade-offs will strengthen the mathematical foundations of physics-guided diffusion method.

**Code**

The code of our article can be found on this website: https://github.com/Prometheus-cotigo/Pde-guide-Diffusion-Model-/tree/main.

**Acknowledgment**

This work was supported by the Hunan Province Natural Science Foundation (2026JJ81495 and 2023JJ70029), and the National Engineering Research Center for High-Speed Railway Construction Technology Open Fund Sponsorship (HSR202211).